\documentclass[11pt]{article}%
\usepackage{amsfonts}
\usepackage{amsmath}
\usepackage{amssymb}
\usepackage{graphicx}%
\providecommand{\U}[1]{\protect\rule{.1in}{.1in}}
\newtheorem{theorem}{Theorem}

\newtheorem{lemma}[theorem]{Lemma}

\newtheorem{remark}[theorem]{Remark}

\newenvironment{proof}[1][Proof]{\noindent\textbf{#1.} }{\ \rule{0.5em}{0.5em}}
\begin{document}

\author{Yu-Chu Lin,\ Chi-Cheung Poon,\ Dong-Ho Tsai\thanks{Research supported by NCTS
and NSC\ of Taiwan under grant number\ 96-2115-M-007-010-MY3.}}
\title{CONTRACTING CONVEX IMMERSED CLOSED\ PLANE CURVES\ WITH\ SLOW\ SPEED OF
CURVATURE\thanks{AMS\ Subject Classifications:\ 35K15, 35K55.}}
\maketitle

\begin{abstract}
We study the contraction of a convex immersed plane curve with speed $\frac
{1}{\alpha}k^{\alpha},\ $where $\alpha\in(0,1]\ $is a constant and show that,
if the blow-up rate of the curvature is of type one, it will converge to a
\textbf{homothetic self-similar solution}.\ We also discuss a special
symmetric\ case of type two blow-up and show that it converges to a
\textbf{translational self-similar solution}. In the case of curve shortening
flow (i.e., when $\alpha=1$), this translational self-similar solution is the
familiar "\textbf{Grim Reaper}"\ (a terminology due to M. Grayson\ \cite{GR}).

\end{abstract}

\section{Introduction.}

Let $\gamma_{0}\ $be\ a convex immersed smooth\ closed plane curve with
\emph{rotation index} (number of times its tangent vector winds around as one
goes along the curve)$\ m\in\mathbb{N},$ parametrized by a smooth immersion
$X_{0}\left(  \varphi\right)  :S^{1}\rightarrow\mathbb{R}^{2}.\ $%
Here\ "convex" means that $\gamma_{0}$ has no inflection points (i.e.,
curvature is positive everywhere). In general, such a\ curve $\gamma_{0}\ $can
have self-intersections (if $m\geq2$). \ 

A family of convex immersed closed curves $X\left(  \varphi,t\right)
:S^{1}\times\lbrack0,T)\rightarrow\mathbb{R}^{2}\ $(with rotation index $m$
and initial data $\gamma_{0}$)\ is said to evolve (contract)\ under the
$k^{\alpha}\ $flow, where $\alpha>0\ $is a constant, if we have
\[
\left(  \bigstar\right)  \ \cdot\cdot\cdot\ \left\{
\begin{array}
[c]{l}%
\dfrac{\partial X}{\partial t}\left(  \varphi,t\right)  =\dfrac{1}{\alpha
}k^{\alpha}\left(  \varphi,t\right)  \mathbf{N}\left(  \varphi,t\right)
,\ \ \ \forall\ \ \left(  \varphi,t\right)  \in S^{1}\times\lbrack0,T)%
\vspace{3mm}%
\\
X\left(  \varphi,0\right)  =X_{0}\left(  \varphi\right)  \in C^{\infty}\left(
S^{1}\right)  ,\ \ \ \varphi\in S^{1},
\end{array}
\right.
\]
where $k\left(  \varphi,t\right)  $ is the curvature of the curve $\gamma
_{t}:=X\left(  \cdot,t\right)  \ $at $\varphi,$ and $\mathbf{N}\left(
\varphi,t\right)  \ $is the unit normal vector of the curve\ $\gamma_{t}%
.\ $Throughout this paper the constant $\alpha$ is assumed to be $0<\alpha
\leq1\ $(in such a case, we shall call $\left(  \bigstar\right)  $
a$\ $\textbf{slow speed\ }contraction).$\ $Here we use the convention that for
convex plane curves\ the curvature$\ k>0$ is positive\ everywhere\ and as for
the direction of the normal $\mathbf{N},\mathbf{\ }$we choose $\mathbf{N}%
=\left(  0,1\right)  $ at a point with minimum $y$-coordinate and extend it
continuously to the whole curve.\ 

When $\alpha=1\ $(i.e., the well-known \emph{curve-shortening flow}),$\ $our
setting is exactly the same as in the interesting paper by Angenent
\cite{ANG}, where the flow $\left(  \bigstar\right)  $ contracts $\gamma_{0}$
with singularity forming in finite time.\ 

Our aim is to study the asymptotic behavior of the contracting flow $\left(
\bigstar\right)  $ with $m\geq2\ $(the immersed case), trying to generalize
results in \cite{ANG} to the case $\alpha\in(0,1].\ $The behavior of the
contracting flow$\ \left(  \bigstar\right)  $\ with \textbf{fast speed}, i.e.,
when$\ 1<\alpha<\infty,$ has been discussed in \cite{PT}. Note that $\alpha
\in(1,\infty)$ corresponds to $p\in\left(  1,2\right)  $ in the equation
$\left(  \clubsuit\right)  $ below.

When $m=1\ $(the embedded case), the initial curve $X_{0}\ $is embedded and
the convergence behavior of the flow $\left(  \bigstar\right)  $ for general
$\alpha>0$ is well understood due to a series of nice papers by Ben
Andrews\ \cite{AN1},\ \cite{AN3}\ and \cite{AN4}. For the information of the
readers,\ we give a brief summary provided by Andrews\footnote{We thank Ben
Andrews for giving us this summary.}:

\begin{theorem}
(\emph{Ben Andrews}\textsf{\ }\cite{AN1},\ \cite{AN3},\ \cite{AN4}%
)\label{thm-ben}\ For $m=1$\ and any $\alpha>0,$ the curve $\gamma_{t}$
contracts to a point in finite time. If $0<\alpha<1/3,$ then for generic
initial data there is no limit of the curves rescaled about the final point
(the isoperimetric ratio approaches infinity); and the exceptional ones where
the isoperimetric ratio remains bounded converge to homothetic solutions,
which have been classified. For $\alpha>1/3,$ the rescaled solutions converge
to circles; and for $\alpha=1/3,$ they converge to ellipses.
\end{theorem}

\begin{remark}
As a consequence of Theorem \ref{thm-ben}, we have the following interesting
elliptic\ result. For $0<\lambda<3\ $(here $\lambda=1/\alpha$)$,$ the only
positive $2\pi$-periodic solution to the equation%
\begin{equation}
w^{\lambda}\left(  x\right)  \left[  w_{xx}\left(  x\right)  +w\left(
x\right)  \right]  =1,\ \ \ x\in S^{1} \label{w-lamda}%
\end{equation}
is $w\left(  x\right)  \equiv1.$ But for $\lambda\geq3,$ we begin to have
nonconstant $2\pi$-periodic\ solutions. For example, when $\lambda=3,$ there
is a family of positive $2\pi$-periodic solutions to equation (\ref{w-lamda})
of the form ($b\in\mathbb{R}\ $is a parameter)%
\begin{equation}
w(x)=\left(  \frac{1}{1+b^{2}}\right)  ^{1/4}\sqrt{1+b^{2}\cos^{2}%
x},\;\;\ x\in\lbrack0,2\pi] \label{wb}%
\end{equation}
where we obtain (\ref{wb}) by computing the curvature (or support
function)\ of an ellipse.\ The function $w\left(  x\right)  \ $of (\ref{wb})
has maximum value occurred at $x=0,\ $with value$\ \left(  1+b^{2}\right)
^{1/4}\geq1$.
\end{remark}

Let $x\in S_{m}^{1}:=\mathbb{R}/2m\pi\mathbb{Z}\ $be the tangent angle of
$\gamma_{t}$ (a\ function defined on $S_{m}^{1}\ $means that it is periodic
with period $2m\pi$).\ In terms of the variable $\left(  x,t\right)  $, it is
known that the curvature quantity $v\left(  x,t\right)  =k^{\alpha}\left(
x,t\right)  \ $of $\gamma_{t}\ $in $\left(  \bigstar\right)  $ will satisfy
the quasilinear parabolic PDE$\ $(a function defined on $S_{m}^{1}\ $means
that it is periodic with period $2m\pi$)$\ $
\[
\left(  \clubsuit\right)  \ \cdot\cdot\cdot\ \left\{
\begin{array}
[c]{l}%
\dfrac{\partial v}{\partial t}=v^{p}\left(  v_{xx}+v\right)  ,\ \ \ p=1+\dfrac
{1}{\alpha}\in\lbrack2,\infty),\ \ \ 0<\alpha\leq1%
\vspace{3mm}%
\\
v\left(  x,0\right)  =v_{0}\left(  x\right)  >0\ \ \ \text{for all\ \ \ }x\in
S_{m}^{1}%
\vspace{3mm}%
\\
v\left(  x,t\right)  =v\left(  x+2m\pi,t\right)  \ \ \ \text{for
all\ \ \ }x\in\mathbb{R},\ \ \ t>0,
\end{array}
\right.
\]
where$\ k_{0}\left(  x\right)  $ is the curvature of $\gamma_{0}$
and$\ v_{0}(x)=k_{0}^{\alpha}\left(  x\right)  ,\ x\in S_{m}^{1}.\ $Moreover,
it is also known that $\left(  \clubsuit\right)  \ $is equivalent to $\left(
\bigstar\right)  \ $(under the assumption that $v_{0}\left(  x\right)  >0$
satisfies the integral condition (\ref{integral-cond}) below).$\ $As
$\gamma_{0}\ $is a closed curve, in $\left(  \bigstar\right)  $\ the initial
data $v_{0}\left(  x\right)  =k_{0}^{\alpha}\left(  x\right)  >0$ in $\left(
\clubsuit\right)  $ must satisfy$\ $the integral condition%
\begin{equation}
\int_{S_{m}^{1}}v_{0}^{1-p}\left(  x\right)  e^{ix}dx=0,\ \ \ e^{ix}=\cos
x+i\sin x \label{integral-cond}%
\end{equation}
where $\int_{S_{m}^{1}}\ $means $\int_{-m\pi}^{m\pi}.\ $Also note that
(\ref{integral-cond})\ is preserved under $\left(  \clubsuit\right)  $, i.e.,
if initially $v_{0}\left(  x\right)  \ $satisfies\ (\ref{integral-cond}), so
does $v\left(  x,t\right)  .$

>From now on we shall focus on $\left(  \clubsuit\right)  \ $with $p\in
\lbrack2,\infty)\ $and $m\geq2\ $with the smooth\ initial data $v_{0}\left(
x\right)  >0\ $in$\ \left(  \clubsuit\right)  \ $satisfying
(\ref{integral-cond}).\ In Lemma \ref{lem3} we shall discuss a result when the
initial function $v_{0}(x)\ $does not satisfy the integral condition
(\ref{integral-cond}). The overall understanding is that when
(\ref{integral-cond}) is satisfied, then we are talking about the geometric
flow $\left(  \bigstar\right)  .$ If not, then one can simply view $\left(
\clubsuit\right)  $\ as a pure analytical problem.\ 

Since equation $\left(  \clubsuit\right)  $\ is parabolic, regularity theory
implies the existence of a unique smooth periodic solution $v\left(
x,t\right)  $ on $S_{m}^{1}\times\lbrack0,T)$ for some $T>0.\ $Each $v\left(
\cdot,t\right)  ,\ t\in\lbrack0,T),$ remains smooth, positive, and periodic
over $\mathbb{R\ }$with period $2m\pi.$ By the equivalence, the flow $\left(
\bigstar\right)  $ also has short time existence of a smooth\ solution.\ Each
$\gamma_{t}$ remains convex, closed, and immersed with rotation index $m$ for
all $t\in\lbrack0,T).\ $

The classical \emph{curve-shortening flow} is when $\alpha=1\ $(or $p=2$); see
Gage-Hamilton \cite{GH} for the embedded case (i.e.,$\ m=1$)$\ $and Angenent
\cite{ANG}, Angenent-Vel\'{a}zquez\ \cite{AV}\ for the immersed case
(i.e.,$\ m\geq2$).\ When $m=1,$ the value $\alpha=1/3$ in Theorem
\ref{thm-ben} corresponds to $p=4$ in $\left(  \clubsuit\right)  $. For more
information on the evolution\ (expansion or contraction) of convex closed
curves in $\mathbb{R}^{2}$, see Andrews \cite{AN1},\ Chou-Zhu\ \cite{CZ}, and
the references therein.\ 

\begin{remark}
If in $\left(  \clubsuit\right)  \ $the constant $\alpha\ $is negative\ (let
$\alpha=-\beta,\ \beta>0$), then the corresponding flow in $\left(
\bigstar\right)  $\ is to expand$\ \gamma_{0}\ $along its outward normal
vector direction with speed $1/\left(  \beta k^{\beta}\right)  $. The
evolution of $v=1/k^{\beta}\ $is\ given by
\begin{equation}
\dfrac{\partial v}{\partial t}=v^{p}(v_{xx}+v),\ \ \ \ \ v\left(  x,t\right)
=v\left(  x+2m\pi,t\right)  ,\ \ \ \ \ p=1+\frac{1}{\alpha}\in\left(
-\infty,1\right)  . \label{v-minus}%
\end{equation}
Finally if one expands $\gamma_{0}$ along its outward normal vector direction
with the exponential speed $\exp\left(  1/k\right)  $, the evolution of
$v=e^{1/k}\ $is%
\begin{equation}
\frac{\partial v}{\partial t}=v\left(  v_{xx}+v\right)  ,\ \ \ \ \ v\left(
x,t\right)  =v\left(  x+2m\pi,t\right)  , \label{v-1}%
\end{equation}
which fills in the gap for $p=1.$\ Hence $p=1$ in $\left(  \clubsuit\right)
\ $separates the contraction case from the expansion case. \ 
\end{remark}

When $m\geq2,\ $the behavior of solutions $v\left(  x,t\right)  \ $to the
equation$\ $%
\begin{equation}
\dfrac{\partial v}{\partial t}=v^{p}\left(  v_{xx}+v\right)
,\ \ \ \ \ v\left(  x,0\right)  =v_{0}\left(  x\right)  >0\in C^{\infty
}\left(  S_{m}^{1}\right)  ,\ \ \ \ \ p\in\left(  -\infty,\infty\right)
\label{v-gen}%
\end{equation}
for $p\in(-\infty,0],$\ $p\in\left(  0,1\right)  ,\ p=1,\ p\in\left(
1,2\right)  ,\ p\in\lbrack2,\infty)$\ are all quite different.\ This means
that equation\ (\ref{v-gen})\ has at least\textbf{ }the following "critical
values".$\ $Each case has its own feature explained below.\ 

\begin{itemize}
\item $p=0.\ $The case of the linear heat equation for the function$\ e^{-t}%
v$, or the case of expanding flow with speed $1/k$.\ It also separates the
sublinear case ($p<0$)\ and the superlinear case\ ($p>0$).

\item $p=1.\ $The case which separates the contraction case from the expansion
case.$\ $In such a case,\ (\ref{integral-cond}) becomes%
\begin{equation}
\int_{S_{m}^{1}}\log v_{0}\left(  x\right)  e^{ix}dx=0. \label{log}%
\end{equation}
Since the behavior of $\log x$ is different\ from $x^{1-p}$ for $p\neq1,$ this
case is quite special.

\item $p=2.\ $The case of the classical \emph{curve-shortening flow}.$\ $It is
the gradient flow of the length functional.\ As we shall see below, for
$p\geq2,\ $(\ref{v-gen}) begins to have type-two blow-up (or type-two
singularity in $\left(  \bigstar\right)  $). Thus $p=2$ separates the type-one
blow-up and the type-two\ blow-up\ \ (see the definition for type-one and
type-two\ blow-up below).
\end{itemize}

\begin{remark}
By Andrews's Theorem \ref{thm-ben}, one can also view $p=4$ as a critical
value although it is for $m=1.$ Also see the discussions before Remark
\ref{rmk1}.\ 
\end{remark}

The behavior of solutions of (\ref{v-gen})\ for $p\geq2$ is most unknown to
us, especially the blow-up rate of a type-two singularity.\ The case $p=1$ is
also complicated.\ Some proofs valid for $p\neq1$ can not be carried over to
the case $p=1.$ To see their differences, we refer to the papers by
Angenent\ \cite{ANG} ($p=2$), Angenent-Vel\'{a}zquez\ \cite{AV}\ ($p=2$%
),\ \cite{PT}\ ($1<p<2$), \cite{T3}\ ($p=1$), \cite{LPT}\ ($0<p<1$)\ and
Urbas\ \cite{U1}\ ($p\leq0$) for details.\ 

Let$\ v_{\min}\left(  t\right)  =\min_{x\in S_{m}^{1}}v\left(  x,t\right)
\ $and similarly for $v_{\max}\left(  t\right)  .\ $For all $p\in\left(
-\infty,\infty\right)  \ $in equation\ (\ref{v-gen}),\ as long as solution
exists,$\ v_{\min}\left(  t\right)  \ $is always increasing due to the maximum
principle.\ By the parabolic regularity theory, it is also known that smooth
solution $v\left(  x,t\right)  $ to equation (\ref{v-gen})\ exists on some
maximal time interval $[0,T_{\max}),$ where $v_{\max}\left(  t\right)
\ $blows up at $T_{\max}\ $($v_{\max}\left(  t\right)  $ will be eventually
increasing for $t$ close to $T_{\max}$).\ For $p\leq0,\ T_{\max}=\infty$\ and
for $p>0,\ T_{\max}<\infty$.

When $T_{\max}<\infty\ $(i.e., when $p>0$), if we let $R\left(  t\right)  $ be
the unique solution to the ODE%
\begin{equation}
\frac{dR}{dt}=R^{p+1}\left(  t\right)  ,\ \ \ \ \ R\left(  T_{\max}\right)
=\infty\label{R-ode}%
\end{equation}
then $R\left(  t\right)  =\left[  p\left(  T_{\max}-t\right)  \right]
^{-1/p}$ and the comparison principle implies that$\ $
\begin{equation}
0<v_{\min}\left(  t\right)  \leq R\left(  t\right)  \leq v_{\max}\left(
t\right)  \ \ \ \text{for all}\ \ \ t\in\lbrack0,T_{\max}). \label{vR}%
\end{equation}
We define the following terminology:\ if there exists a constant $C,$
independent\ of time, such that%
\begin{equation}
0<v_{\max}\left(  t\right)  \leq CR\left(  t\right)  \ \ \ \text{for
all}\ \ \ t\in\lbrack0,T_{\max}) \label{t-one}%
\end{equation}
then we say $v\left(  x,t\right)  $ has \textbf{type-one} blow-up. If not,
i.e.,$\ $if$\ v_{\max}\left(  t\right)  /R\left(  t\right)  \ $is not
bounded\ on$\ t\in\lbrack0,T_{\max}),\ $then we say $v\left(  x,t\right)  $
has \textbf{type-two} blow-up. A type-two blow-up is clearly much more complicated.\ 

It has been shown in p. 158 of \cite{LPT} that for $p\in\left(  0,2\right)  $
there is \textbf{no} type-two blow-up for any $m\in\mathbb{N}\ $(include
$m=1$)\ and any positive\ initial data $v_{0}\left(  x\right)  \in C^{\infty
}\left(  S_{m}^{1}\right)  \ $(no matter it satisfies (\ref{integral-cond}) or
not), i.e., all blow-ups are of type-one in $\left(  \clubsuit\right)  \ $for
$p\in\left(  0,2\right)  \ $and for any $v_{0}\left(  x\right)  >0$%
.$\ $However,\ for $p\in\left(  0,2\right)  ,$ the limit function $w\left(
x\right)  \ $($w\left(  x\right)  $ is the limit of the rescaled
solution\ $v\left(  x,t\right)  /R\left(  t\right)  $)$\ $may be either
$w\left(  x\right)  >0\ $everywhere on\ $S_{m}^{1}\ $(call it
\textbf{nondegenerate}) or $w\left(  x\right)  \geq0\ $but with $w\left(
x\right)  =0$\textsf{ }somewhere on\ $S_{m}^{1}\ $(call it \textbf{degenerate}%
). If $w\left(  x\right)  $ is nondegenerate, it gives rise to a
\textbf{self-similar (homothetic)\ solution}.\ If $w\left(  x\right)  $ is
degenerate, its behavior (regularity)\ for $p\in\left(  0,1\right)
,\ p=1,\ $and$\ p\in\left(  1,2\right)  $ are all different.\ 

In Angenent \cite{ANG} (the case $p=2$), he employed an elegant unstable
manifold analysis of a shrinking spiral (travelling wave solution) and used it
to prove a\textsf{ }Harnack-type estimate,\ i.e.,\textsf{\ }Lemma 7.1. of
\cite{ANG}. This is the key estimate to ensure convergence to a positive
self-similar solution under type-one blow-up\ (see p. 605,\ Theorem A of his paper).\ 

We shall see that his proof can be carried to the case $p\in\lbrack
2,\infty),\ $assuming that we have type-one blow-up.\ Hence we can obtain
convergence of equation $\left(  \clubsuit\right)  $ (or the flow $\left(
\bigstar\right)  $)\ to a self-similar solution$\ w\left(  x\right)  $ under
type-one blow-up. Unlike the case for $p\in\left(  0,2\right)  ,\ $now
$w\left(  x\right)  \ $is positive everywhere on $S_{m}^{1}\ $(for
$p\in\lbrack2,\infty)$) and is an entire periodic solution to the
corresponding steady state\ ODE. See Theorem \ref{thmA} in Section
\ref{type-one}.

In summary, the above says that for $p\in\lbrack2,\infty)\ $there is either
type-one blow-up\ or type-two blow-up.\ Moreover, if we have type-one blow-up,
then the limit function $w\left(  x\right)  $ is always \emph{nondegenerate},
i.e., $w\left(  x\right)  >0$ everywhere.\ 

As for Theorem\ B of \cite{ANG},\ it is rather straightforward to generalize
it to the case $p\geq2.\ $See Theorem\ \ref{thmB}.

Finally we also discuss a special \textbf{symmetric}\ case of type-two blow-up
and obtain a convergence to the cosine function (see Theorem \ref{thmC}\ and
Theorem \ref{thmC-1}). In flow\ $\left(  \bigstar\right)  $, this convergence
gives rise to a \textbf{translational self-similar solution}. When $\alpha
=1$\ (curve-shortening flow), this translational self-similar solution is the
Grayson's "\textbf{grim reaper}", i.e., the graph of $y=-\log\cos x,$
$x\in\left(  -\pi/2,\pi/2\right)  .\ $One can view Theorem \ref{thmC}\ and
Theorem \ref{thmC-1}\ as partial generalizations of Theorem\ C of \cite{ANG}
to the case $p\geq2.\ $

In conclusion, we can generalize Theorem A, Theorem\ B,\ and part of Theorem C
in p. 605 of Angenent \cite{ANG} to the case $p\geq2.\ $

To end this introductory section, we point out that solutions of (\ref{v-gen})
for the sublinear case $p<0$ are well-behaved as it is bounded by the
following super-sub solutions
\begin{equation}
\left[  \left(  \min_{x\in S_{m}^{1}}v_{0}\left(  x\right)  \right)
^{-p}-pt\right]  ^{-1/p}\leq v\left(  x,t\right)  \leq\left[  \left(
\max_{x\in S_{m}^{1}}v_{0}\left(  x\right)  \right)  ^{-p}-pt\right]  ^{-1/p}
\label{minmax}%
\end{equation}
for all$\ t\in\lbrack0,\infty).\ $In particular we have%
\[
1\leq\frac{v_{\max}\left(  t\right)  }{v_{\min}\left(  t\right)  }\leq
U\left(  t\right)  :=\frac{\left[  \left(  \max_{x\in S_{m}^{1}}v_{0}\left(
x\right)  \right)  ^{-p}-pt\right]  ^{-1/p}}{\left[  \left(  \min_{x\in
S_{m}^{1}}v_{0}\left(  x\right)  \right)  ^{-p}-pt\right]  ^{-1/p}}%
,\ \ \ t\in\lbrack0,\infty)
\]
where$\ U\left(  t\right)  $ is a \textbf{decreasing} function on$\ [0,\infty
)\ $with $\lim_{t\rightarrow\infty}U\left(  t\right)  =1.\ $In fact, the
quantity $v_{\max}\left(  t\right)  /v_{\min}\left(  t\right)  $ also
decreases to $1$\ as $t\rightarrow\infty$. As a consequence, by regularity
theory, the rescaled solution$\ u\left(  x,t\right)  :=v\left(  x,t\right)
/R\left(  t\right)  $ will converge as $t\rightarrow\infty\ $to the constant
function $1$ in $C^{k}\left(  S_{m}^{1}\right)  \ $for any $k\in\mathbb{N}.$
Here $R\left(  t\right)  $ can be \emph{any}\ solution to the ODE
$dR/dt=R^{1+p}\ $($p<0$)$\ $with $R\left(  0\right)  >0.$ The geometric
meaning is that when $\alpha\in\left(  -1,0\right)  ,$ the expanding flow
$\left(  \bigstar\right)  $ converges (after rescaling)\ to the $m$-fold cover
of $S^{1}\ $in any $C^{k}$-norm.\ See Urbas \cite{U1} also.\ 

\section{Some basic estimates.}

>From now on we assume $p\geq2\ $and $m\geq2,\ $with the smooth\ initial data
$v_{0}\left(  x\right)  >0\ $in$\ \left(  \clubsuit\right)  \ $satisfying
(\ref{integral-cond}).\ For convenience, denote the maximal space-time domain
$S_{m}^{1}\times\lbrack0,T_{\max})\ $by $\Omega_{m}$. In below, if the proof
of a lemma is omitted, then it is either straightforward or similar to those
established in \cite{GH} or \cite{ANG} for the case $p=2$. Hence we will not
repeat it.\ 

\begin{lemma}
\label{lem1}(\emph{gradient estimate in integral form}) There exists a
constant $C$ depending only on $v_{0}$ such that%
\begin{equation}
\int_{S_{m}^{1}}v_{x}^{2}\left(  x,t\right)  dx\leq\int_{S_{m}^{1}}%
v^{2}\left(  x,t\right)  dx+C \label{grad-int}%
\end{equation}
for all $t\in\lbrack0,T_{\max}),\ $where $v_{x}^{2}$ means $\left(  \partial
v/\partial x\right)  ^{2}.\ $In particular, for $\varepsilon>0$
sufficiently\ small,\ there exists\ a number\ $\delta>0,\ $depending only on
$\varepsilon$, such that%
\begin{equation}
\left(  1-\varepsilon\right)  v_{\max}\left(  t\right)  \leq v\left(
x,t\right)  +\sqrt{2m\pi C} \label{local-Har}%
\end{equation}
for all\ $x\in\left(  x_{t}-\delta^{2},x_{t}+\delta^{2}\right)  \ $and all
$t\in\lbrack0,T_{\max}),\ $where$\ v\left(  x_{t},t\right)  =v_{\max}\left(
t\right)  .\ $
\end{lemma}

\begin{lemma}
\label{lem2}(\emph{gradient estimate})\ For solution $v\left(  x,t\right)
\ $to equation $\left(  \clubsuit\right)  $\ on $\Omega_{m}\ $we have, at each
point $(x,t),$\ either
\begin{equation}
\left(  v_{xx}+v\right)  \left(  x,t\right)  >0 \label{Ang2}%
\end{equation}
or
\begin{equation}
v_{x}^{2}\left(  x,t\right)  +v^{2}\left(  x,t\right)  \leq\max_{x\in
S_{m}^{1}}\left[  \left(  v_{0}\right)  _{x}^{2}\left(  x\right)  +v_{0}%
^{2}\left(  x\right)  \right]  :=\sigma. \label{Ang3}%
\end{equation}
In particular we have
\begin{equation}
\left\vert v_{x}\left(  x,t\right)  \right\vert \leq\max\left\{
\lambda,\ v_{\max}\left(  t\right)  \right\}  \label{vvv}%
\end{equation}
for all $\left(  x,t\right)  \in\Omega_{m},$ where $\lambda>0$ is a constant
depending only on $v_{0}$.\ As a consequence we also know that $v_{\max
}\left(  t\right)  $ is eventually increasing for $t$ close to $T_{\max}.\ $
\end{lemma}

\begin{lemma}
\label{lem2-1}(\emph{behavior near maximum point})\ Let $v_{\max}\left(
t\right)  =v\left(  x_{t},t\right)  \ $for some $x_{t}\in S_{m}^{1}$.$\ $If at
any time $t\in\lbrack0,T_{\max})$ we have $v_{\max}(t)>\sigma,$\ where
$\sigma$ is from (\ref{Ang3}),\ then
\begin{equation}
v\left(  x,t\right)  >v_{\max}\left(  t\right)  \cos\left(  x-x_{t}\right)
\label{vvcos}%
\end{equation}
for all $x$\ with $0<\left\vert x-x_{t}\right\vert <\arccos\left(
\sigma/v_{\max}\left(  t\right)  \right)  .$
\end{lemma}

With the help of the above basic estimates, we can generalize Theorem B of
\cite{ANG} to the case $p\geq2:$

\begin{theorem}
\label{thmB}(\emph{rough upper bound of}\textsf{ }$v_{\max}\left(  t\right)
$)\ If $v_{\max}\left(  t\right)  $ blows up at time $T_{\max},$ then there
holds the following
\begin{equation}
\lim_{t\rightarrow T_{\max}}\left(  T_{\max}-t\right)  v_{\max}\left(
t\right)  =0. \label{rough-est}%
\end{equation}

\end{theorem}

\begin{remark}
For now, by (\ref{vR})\ and (\ref{rough-est}) we have the rough estimate%
\[
\left[  p\left(  T_{\max}-t\right)  \right]  ^{-1/p}\leq v_{\max}\left(
t\right)  \leq\frac{C}{T_{\max}-t}\ \ \ \text{for all\ \ \ }t\in
\lbrack0,T_{\max}),
\]
where$\ p\geq2$ and $C$ is some constant independent\ of time.\ 
\end{remark}

%

\proof
Let
\[
I\left(  t\right)  =\int_{S_{m}^{1}}v^{1-p}\left(  x,t\right)
dx>0,\ \ \ p\geq2.
\]
By (\ref{local-Har})\ in Lemma \ref{lem1}\ we have for $t$ close to $T_{\max}$
the estimate
\[
\int_{S_{m}^{1}}v\left(  x,t\right)  dx\geq cv_{\max}\left(  t\right)
\]
where $c>0$ is a constant independent\ of time.\ Hence there is a time
$t_{\ast}$ close to $T_{\max}\ $such that
\[
-I^{\prime}\left(  t\right)  =\left(  p-1\right)  \int_{S_{m}^{1}}v\left(
x,t\right)  dx\geq\left(  p-1\right)  cv_{\max}\left(  t\right)
>0\ \ \ \text{for all\ \ \ }t\in\lbrack t_{\ast},T_{\max})
\]
By integration of $I^{\prime}\left(  t\right)  $ on the interval $[t_{\ast
},T_{\max}),$ we obtain%
\[
\left(  p-1\right)  c\int_{t_{\ast}}^{T_{\max}}v_{\max}\left(  t\right)
dt\leq I\left(  t_{\ast}\right)  -\left(  \lim_{t\rightarrow T_{\max}}I\left(
t\right)  \right)  \leq I\left(  t_{\ast}\right)  <\infty
\]
and so the integral$\ \int_{0}^{T_{\max}}v_{\max}\left(  t\right)  dt\ $is
finite. Since by Lemma \ref{lem2} $v_{\max}\left(  t\right)  $ is eventually
increasing, we may also assume that $v_{\max}\left(  t\right)  $ is increasing
on $[t_{\ast},T_{\max})$ and conclude
\begin{equation}
\left(  T_{\max}-t\right)  v_{\max}\left(  t\right)  \leq\int_{t}^{T_{\max}%
}v_{\max}\left(  s\right)  ds\ \ \ \text{for all\ \ \ }t\in\lbrack t_{\ast
},T_{\max}). \label{ttt}%
\end{equation}
Letting $t\rightarrow T_{\max},$ the right hand side of (\ref{ttt}) converges
to zero and the proof is done.$%
\hfill
\square$

\section{Type-one blow-up implies $C^{\infty}$ convergence. \label{type-one}}

Throughout this section we assume the solution $v\left(  x,t\right)  $ to
equation $\left(  \clubsuit\right)  $ has type-one blow-up. We shall consider
its asymptotic\ behavior by the obvious rescaling $u\left(  x,t\right)
:=v\left(  x,t\right)  /R\left(  t\right)  ,\ $where\ $R\left(  t\right)
\ $is from\ (\ref{R-ode}),\ and let $\tau\in\lbrack0,\infty)\ $be the new time
given by the relation$\ t=T_{\max}\left(  1-e^{-p\tau}\right)  ,\ t\in
\lbrack0,T_{\max}),\ $which is motivated by the requirement $d\tau
/dt=R^{p}\left(  t\right)  ,\ $then the function
\begin{equation}
u\left(  x,\tau\right)  =p^{1/p}T_{\max}^{1/p}e^{-\tau}v\left(  x,T_{\max
}\left(  1-e^{-p\tau}\right)  \right)  >0,\ \ \ x\in S_{m}^{1},\ \ \ \tau
\in\lbrack0,\infty) \label{rescale}%
\end{equation}
\ will be a positive, \textbf{bounded}, solution of the rescaled equation%
\begin{equation}
\left\{
\begin{array}
[c]{l}%
\dfrac{\partial u}{\partial\tau}=u^{p}\left(  u_{xx}+u-u^{1-p}\right)
,\ \ \ p\geq2%
\vspace{3mm}%
\\
u\left(  x,\tau\right)  =u(x+2m\pi,\tau)
\end{array}
\right.  \label{dudt}%
\end{equation}
for all$\ \left(  x,\tau\right)  \in S_{m}^{1}\times\lbrack0,\infty),\ $with
$u\left(  x,0\right)  =u_{0}\left(  x\right)  :=p^{1/p}T_{\max}^{1/p}%
v_{0}\left(  x\right)  >0.$ Moreover, we have%
\begin{equation}
0<u_{\min}\left(  \tau\right)  \leq1\leq u_{\max}\left(  \tau\right)
\label{UU}%
\end{equation}
for\ all$\ \tau\in\lbrack0,\infty)\ $due to (\ref{vR}). By (\ref{vvv}) we also
have the uniform gradient estimate%
\begin{equation}
\left\vert u_{x}\left(  x,\tau\right)  \right\vert \leq C\ \ \ \text{for
all\ \ \ }\left(  x,\tau\right)  \in S_{m}^{1}\times\lbrack0,\infty)
\label{grad-est}%
\end{equation}
where $C$\ is a constant depending only on $v_{0}.$

We shall generalize Angenent's Lemma 7.1 in\ \cite{ANG} to the following:

\begin{theorem}
\label{thm-type-one}(\emph{gradient estimate for type-one blow-up})\ Let
$v\left(  x,t\right)  $ be a \emph{type-one solution} to equation $\left(
\clubsuit\right)  $\ with $p\geq2$. Then the rescaled
bounded\ positive\ function $u\left(  x,\tau\right)  $ of (\ref{rescale})
satisfies \
\begin{equation}
\left\vert u_{x}\left(  x,\tau\right)  \right\vert \leq\lambda u\left(
x,\tau\right)  \ \ \ \text{for all\ \ \ }\left(  x,\tau\right)  \in S_{m}%
^{1}\times\lbrack0,\infty) \label{KK}%
\end{equation}
where $\lambda$ is a constant depending only on$\ u_{0}$.
\end{theorem}

\begin{remark}
Theorem \ref{thm-type-one}$\ $fails for $p\in\left(  0,2\right)  .\ $
\end{remark}

\begin{remark}
By (\ref{KK})\ we have the estimate%
\begin{equation}
1\leq u_{\max}\left(  \tau\right)  \leq e^{2\lambda m\pi}u_{\min}\left(
\tau\right)  \label{KK1}%
\end{equation}
for all $\tau\in\lbrack0,\infty)\ $and hence $u_{\min}\left(  \tau\right)  $
has a positive lower bound for $\tau\in\lbrack0,\infty)$ and equation
(\ref{dudt}) is uniformly parabolic\ on $S_{m}^{1}\times\lbrack0,\infty)$.
\end{remark}

\begin{remark}
>From the proof we see that Theorem \ref{thm-type-one} remains valid even the
initial condition $v_{0}\left(  x\right)  $ does not satisfy the integral
condition (\ref{integral-cond}).\ This observation is important and will be
used in Lemma \ref{lem3} below.
\end{remark}

\subsection{Angenent's method of shrinking spirals.}

Since Theorem \ref{thm-type-one}\ is valid for $p=2$, we assume $p>2.\ $Our
method of proof is similar to that originally used by Angenent in
\cite{ANG}.\ At the same time we also provide some additional details and see
why we need the condition $p>2.\ $Consider a special solution
(\emph{travelling wave solution}) of the form $U\left(  x,\tau\right)
=U\left(  x-c\tau\right)  ,\ c>0\ $($c\ $is a constant), to the equation%
\begin{equation}
\dfrac{\partial u}{\partial\tau}=u^{p}\left(  u_{xx}+u-u^{1-p}\right)  .
\label{eq-ss}%
\end{equation}
A positive function $U\left(  \xi\right)  \ $over some interval $I$ will
generate a solution if and only if
\begin{equation}
U^{p}\left(  \xi\right)  U^{\prime\prime}\left(  \xi\right)  +U^{p+1}\left(
\xi\right)  -U\left(  \xi\right)  +cU^{\prime}\left(  \xi\right)
=0\ \ \ \text{for all\ \ \ }\xi\in I. \label{Q1}%
\end{equation}
For such a$\ U\left(  \xi\right)  >0$ satisfying equation (\ref{Q1}) on $I$ we
have
\begin{equation}
\frac{d}{d\xi}\left[  \left(  U^{\prime}\left(  \xi\right)  \right)
^{2}+U^{2}\left(  \xi\right)  -\frac{2}{2-p}U^{2-p}\left(  \xi\right)
\right]  =-2c\frac{\left(  U^{\prime}\left(  \xi\right)  \right)  ^{2}}%
{U^{p}\left(  \xi\right)  }\leq0 \label{Lyapunov}%
\end{equation}
and so the function$\ E\left(  \xi\right)  $ given by
\begin{equation}
E\left(  \xi\right)  :=\left(  U^{\prime}\left(  \xi\right)  \right)
^{2}+U^{2}\left(  \xi\right)  -\frac{2}{2-p}U^{2-p}\left(  \xi\right)
,\ \ \ \xi\in I,\ \ \ p>2, \label{E}%
\end{equation}
is decreasing in $\xi\in I\ $if $c\neq0.\ $For$\ c\neq0$ the only periodic
solution for (\ref{Q1})\ is the constant$\ U\left(  \xi\right)  \equiv1.$

For $c=0,$ $E\left(  \xi\right)  $ is a positive constant independent\ of
$\xi\in I.\ $It is obvious that any positive\ solution $U\left(  \xi\right)  $
satisfying the equation $E\left(  \xi\right)  =const.>0$ can not become too
small over its domain\ since$\ -2\left(  2-p\right)  ^{-1}U^{2-p}\left(
\xi\right)  \rightarrow\infty$ as$\ U\left(  \xi\right)  \rightarrow0^{+}%
.\ $Thus any solution $U\left(  \xi\right)  \ $to the ODE $U^{p}%
U^{\prime\prime}+U^{p+1}-U=0\ $is a\ positive periodic function on$\ \xi
\in\left(  -\infty,\infty\right)  \ $(here we need the condition
$p>2$)\ satisfying
\begin{equation}
\left(  U^{\prime}\left(  \xi\right)  \right)  ^{2}+U^{2}\left(  \xi\right)
-\frac{2}{2-p}U^{2-p}\left(  \xi\right)  =b^{2}-\frac{2}{2-p}b^{2-p}%
=a^{2}-\frac{2}{2-p}a^{2-p},\ \ \ \xi\in\left(  -\infty,\infty\right)  ,
\label{ODE}%
\end{equation}
where $b\geq1$ ($a\leq1$)\ is the maximum (minimum)\ value of$\ U\left(
\xi\right)  $ over$\ \left(  -\infty,\infty\right)  .$

Similar to Theorem 6.1 of \cite{ANG}, we claim the following:

\begin{theorem}
\label{thm-type-one-1}Assume $p>2.\ $For any small $c>0,$ there is a unique
positive\ solution $U_{c}$ $\in C^{\infty}\left(  (-\infty,0]\right)  \ $of
(\ref{Q1})$\ $with $\lim_{\xi\rightarrow-\infty}U_{c}\left(  \xi\right)
=0\ $and\ the following properties:
\begin{equation}
\left\{
\begin{array}
[c]{l}%
\text{(i).}\ \ U_{c}^{\prime}\left(  \xi\right)  >0\ \ \ \text{for \ \ }\xi
\in(-\infty,0),%
\vspace{3mm}%
\\
\text{(ii).\ \ }U_{c}^{\prime}\left(  0\right)  =0%
\vspace{3mm}%
\\
\text{(iii).\ \ }U_{c}^{\prime}\left(  \xi\right)  \leq\lambda_{c}U\left(
\xi\right)  \ \ \ \text{for \ \ }\xi\in(-\infty,0],
\end{array}
\right.  \label{Q2}%
\end{equation}
where $\lambda_{c}>0\ $is a large constant depending on $c\ $(and
$p\ $also).$\ $As a function of $c>0,$ $U_{c}\left(  0\right)  $ is strictly
decreasing and given any $\delta>0$ and $A>0,$ one can choose $c=c\left(
\delta,A\right)  >0\ $so small that%
\begin{equation}
U_{c}\left(  0\right)  >\delta^{-1}\ \ \ \ \ \ \ \text{and\ \ \ \ \ \ \ }%
U_{c}^{\prime}\left(  \xi\right)  >A\ \ \text{\ whenever\ \ \ }\delta\leq
U_{c}\left(  \xi\right)  \leq\delta^{-1}. \label{Q3}%
\end{equation}

\end{theorem}%

\proof
It will be convenient to look at $H=U^{p}>0$ instead of $U$ itself. We have
$U=H^{1/p}$ and (\ref{Q1}) is equivalent\ to%
\[
HH_{\xi\xi}=-cH_{\xi}+\frac{p-1}{p}\left(  H_{\xi}\right)  ^{2}-pH^{2}%
+pH,\ \ \ H\left(  \xi\right)  =U^{p}\left(  \xi\right)  ,
\]
which can be written as the first order system (let $G=H_{\xi}$)%
\begin{equation}
\left\{
\begin{array}
[c]{l}%
HH_{\xi}=HG%
\vspace{3mm}%
\\
HG_{\xi}=-pH^{2}+pH-cG+\frac{p-1}{p}G^{2}.
\end{array}
\right.  \label{ODE-sys}%
\end{equation}
Thus \emph{up to a reparametrization} (since there is a factor $H$ in front of
$H_{\xi}\ $and $G_{\xi}$), the positive solutions of (\ref{Q1}) are in
one-to-one correspondence with the orbits of the vector field%
\begin{equation}
X_{c}\left(  H,G\right)  =\left(  HG,\ -pH^{2}+pH-cG+\frac{p-1}{p}%
G^{2}\right)  \label{X}%
\end{equation}
lying on the region$\ R^{+}=\left\{  \left(  H,G\right)  :H>0\right\}  .\ $We
shall analyze the phase portrait of the vector field $X_{c}\left(  H,G\right)
\ $in $R^{+}.$ Note that $X_{c}$ has three zeros $\left(  0,0\right)  ,$
$\left(  1,0\right)  ,$ and $\left(  0,pc/\left(  p-1\right)  \right)  \ $in
$R^{+}.\ $For a given small $c>0,\ $our aim is to look at certain special
solution $\left(  H\left(  \xi\right)  ,G\left(  \xi\right)  \right)  $ of the
system (\ref{ODE-sys}) with $H\left(  \xi\right)  >0\ $everywhere.\ 

If we compute the linearization of $X_{c}\left(  H,G\right)  $ at these
equilibrium points, we obtain the three\ matrices
\[
M=%
\begin{pmatrix}
G & H\\
p-2pH & -c+\frac{2\left(  p-1\right)  }{p}G
\end{pmatrix}
=%
\begin{pmatrix}
0 & 0\\
p & -c
\end{pmatrix}
,\ \ \
\begin{pmatrix}
0 & 1\\
-p & -c
\end{pmatrix}
,\ \ \
\begin{pmatrix}
\frac{pc}{p-1} & 0\\
p & c
\end{pmatrix}
\]
at $\left(  0,0\right)  ,$ $\left(  1,0\right)  \ $and $\left(  0,pc/\left(
p-1\right)  \right)  $ respectively.\ The eigenvalues of them are given
respectively by
\[
\lambda=0,\ -c;\ \ \ \ \ \lambda=\frac{-c\pm\sqrt{c^{2}-4p}}{2}%
;\ \ \ \ \ \lambda=\frac{pc}{p-1},\ c
\]
where $c>0$ is a small constant to be chosen later\ on. Therefore, $\left(
0,0\right)  $ is a degenerate zero of $X_{c}$ (it is not a hyperbolic fixed
point of $X_{c}$),$\ \left(  1,0\right)  $ is a spiraling sink of $X_{c}$ if
$c>0$ is small with $c^{2}-4p<0,$ and $\left(  0,pc/\left(  p-1\right)
\right)  $ is a source of $X_{c}.$

By definition, the unstable set $W^{u}\left(  O\right)  $ of the origin
$O=\left(  0,0\right)  \ $consists of all orbits of $X_{c},$ which tend to $O$
as $\xi\rightarrow-\infty.$ As the origin is degenerate, one needs to analyze
further to know what$\ W^{u}\left(  O\right)  $ looks like.\ Note that if
$U_{c}\left(  \xi\right)  $ is the solution satisfying Theorem
\ref{thm-type-one-1}, then
\[
\left(  H_{c}\left(  \xi\right)  ,G_{c}\left(  \xi\right)  \right)  =\left(
U_{c}^{p}\left(  \xi\right)  ,\ pU_{c}^{p-1}\left(  \xi\right)  U_{c}^{\prime
}\left(  \xi\right)  \right)
\]
parametrizes a trajectory of $X_{c}$ in the unstable set $W^{u}\left(
O\right)  \ $of the origin.\ Thus one needs to look at $W^{u}\left(  O\right)
.\ $

\ \ \ \ 

\underline{\textbf{Existence of a trajectory in }$W^{u}\left(  O\right)  .$%
}\ \ Given a constant $\lambda>0$ and let $l_{\lambda}$ be the half line
$G=\lambda H,$ $H\geq0.$ The half line $l_{\lambda}\ $has upward normal
$\left(  -\lambda,1\right)  $ and along it we have
\[
\left\langle X_{c},\left(  -\lambda,1\right)  \right\rangle =\left\{
p-c\lambda-\left(  p+\frac{\lambda^{2}}{p}\right)  H\right\}  H,\ \ \ \lambda
>0,\ \ \ H\geq0.
\]
Now choose two positive$\ \lambda_{1},\ \lambda_{2}$ such that $\lambda
_{1}<p/c\ $and$\ \lambda_{2}=p/c.\ $Put$\ h^{\ast}=\left(  p-c\lambda
_{1}\right)  /\left(  p+\lambda_{1}^{2}/p\right)  >0,$ and define the points
\[
A=\left(  h^{\ast},\lambda_{1}h^{\ast}\right)  ,\ \ \ B=\left(  h^{\ast
},\lambda_{2}h^{\ast}\right)  ,\ \ \ O=\left(  0,0\right)  .
\]
Along the segment $OA\ $with $0\leq H\leq h^{\ast},$ we have $\left\langle
X_{c},\left(  -\lambda_{1},1\right)  \right\rangle \geq0\ $($=0\ $only at
$H=0\ $or$\ H=h^{\ast}$) and so the vector field $X_{c},$ when restricted to
$OA,$ is pointing toward the upper half of segment $OA.\ $Similarly, along the
segment $OB$ we have$\ \left\langle X_{c},\left(  -\lambda_{2},1\right)
\right\rangle \leq0\ $($=0\ $only at $H=0$).\ Hence the vector field $X_{c},$
when restricted to $OB,$ is pointing toward the lower half of the segment.

Finally, along the segment $AB,$ we have$\ \left\langle X_{c},\left(
1,0\right)  \right\rangle =h^{\ast}G>0$ for all $G$ with$\ \lambda_{1}h^{\ast
}<G\ <\lambda_{2}h^{\ast}.\ $We conclude that the\ trajectories\ of$\ X_{c}%
\ $enter\ the\ triangle$\ OAB\ $through\ the\ sides$\ OA\ $and$\ OB,\ $and
they leave\ $OAB\ $through the vertical side $AB.$

For any point $\left(  H,G\right)  $ inside the triangle $OAB,$ it has the
form $G=\lambda H,$ for some $\lambda_{1}\leq\lambda\leq\lambda_{2}=p/c$ and
hence
\begin{align}
G_{\xi}  &  =\frac{1}{H}\left(  -pH^{2}+pH-c\lambda H+\frac{p-1}{p}\lambda
^{2}H^{2}\right) \nonumber\\
&  =\left(  p-c\lambda\right)  +\left(  \frac{p-1}{p}\lambda^{2}-p\right)
H\geq\left(  \frac{p-1}{p}\lambda_{1}^{2}-p\right)  H. \label{G}%
\end{align}
We may choose $c$ small enough and $\lambda_{1}<p/c$ larger than $p/\sqrt
{p-1}$ such that $p/\sqrt{p-1}<\lambda_{1}\leq\lambda\leq\lambda_{2}=p/c\ $and
conclude that
\[
G_{\xi}\geq\left(  \frac{p-1}{p}\lambda_{1}^{2}-p\right)  H>0
\]
for all $\lambda_{1}\leq\lambda\leq\lambda_{2}$ and $0<H\leq h^{\ast}=\left(
p-c\lambda_{1}\right)  /\left(  p+\lambda_{1}^{2}/p\right)  .$ Therefore
$G_{\xi}>0$ in the interior of the triangle $OAB$ and it follows from the
\emph{Wa\.{z}ewski's Principle} that at least one of the trajectories through
$AB$ tends to the origin as $\xi\rightarrow-\infty.$

In conclusion, we see that as long as $c>0$ is small enough (depending only on
$p$), there exists a trajectory in the unstable set\textbf{ }$W^{u}\left(
O\right)  .$

\ \ \ \ \ 

\underline{\textbf{Uniqueness of the trajectory in }$W^{u}\left(  O\right)
.$}\ \ Express the second equation of (\ref{ODE-sys}) as%
\[
HG_{\xi}=-p\left(  H-\frac{1}{2}\right)  ^{2}+\frac{p-1}{p}\left(  G-\frac
{pc}{2\left(  p-1\right)  }\right)  ^{2}+\frac{p}{4}\left(  1-\frac{c^{2}%
}{p-1}\right)  .
\]
We see that if $c<\sqrt{p-1},$ then the set $\left\{  \left(  H,G\right)
:H>0,\ G_{\xi}=0\right\}  $ is the part of the hyperbola
\[
\Gamma:\ p\left(  H-\frac{1}{2}\right)  ^{2}-\frac{p-1}{p}\left(  G-\frac
{pc}{2\left(  p-1\right)  }\right)  ^{2}=\frac{p}{4}\left(  1-\frac{c^{2}%
}{p-1}\right)
\]
lying in $R^{+}=\left\{  \left(  H,G\right)  :H>0\right\}  .\ $Here $\Gamma$
is centered at $\left(  1/2,\ pc/\left(  2\left(  p-1\right)  \right)
\right)  $ and passes through the three equilibrium points $\left(
0,0\right)  ,$ $\left(  0,pc/\left(  p-1\right)  \right)  \ $and $\left(
1,0\right)  .$

Consider the region $\Omega$ enclosed by the segment $\left\{  \left(
0,G\right)  :0\leq G\leq pc/\left(  p-1\right)  \right\}  $ and the left
branch of $\Gamma.\ $For each\ point$\ \left(  H,G\right)  \ $inside
$\Omega,\ $the vector field $X_{c}\left(  H,G\right)  =\left(  HH_{\xi
},HG_{\xi}\right)  \ $satisfies $HH_{\xi}=HG>0\ $and $HG_{\xi}<0.$ This
implies that any trajectory $\left(  H,G\right)  \ $in\ $W^{u}\left(
O\right)  $ will not pass through the region $\Omega$ and thus there exists a
large constant $\lambda_{c}\ $(say$\ \lambda_{c}>p/c,\ $where $p/c$ is the
slope of $\Gamma\ $at the origin) such that any any trajectory $\left(
H\left(  \xi\right)  ,G\left(  \xi\right)  \right)  $ in $W^{u}\left(
O\right)  $ satisfies $G\leq\lambda_{c}H,\ $as $\xi\rightarrow-\infty.$

Let $\left(  H_{1},G_{1}\right)  $ and $\left(  H_{2},G_{2}\right)  \ $be two
different orbits in $W^{u}\left(  O\right)  .\ $From the above observation,
near the origin they can be represented as the graphs $G_{1}=g_{1}\left(
H\right)  ,\ G_{2}=g_{2}\left(  H\right)  ,$ where the $g_{i}$ are solutions
of the equation
\[
g^{\prime}\left(  H\right)  =\frac{HG_{\xi}}{HH_{\xi}}=\frac{-pH^{2}%
+pH-cG+\frac{p-1}{p}G^{2}}{HG}=\frac{\left(  p-1\right)  g\left(  H\right)
-pc}{pH}+\frac{p-pH}{g\left(  H\right)  }.
\]
Orbits cannot intersect, so we may assume that $g_{1}\left(  H\right)
<g_{2}\left(  H\right)  .$Their difference $w\left(  H\right)  =g_{2}\left(
H\right)  -g_{1}\left(  H\right)  >0$ satisfies the equation
\[
w^{\prime}\left(  H\right)  =\left(  \frac{p-1}{pH}-\frac{p-pH}{g_{2}\left(
H\right)  g_{1}\left(  H\right)  }\right)  w,\ \ \ w>0.
\]
Now for $H>0$ sufficiently small, we have $g_{1}\left(  H\right)  \leq\lambda
H,$ $g_{2}\left(  H\right)  \leq\lambda H$ for some $\lambda,$ and so
\[
w^{\prime}\left(  H\right)  \leq\left(  \frac{p-1}{pH}-\frac{p-pH}%
{g_{2}\left(  H\right)  g_{1}\left(  H\right)  }\right)  w=\frac{w}%
{p\lambda^{2}H^{2}}\left[  \lambda^{2}\left(  p-1\right)  H+p^{2}%
H-p^{2}\right]  <0
\]
for all\ sufficiently small $H>0,\ $which means that $w\left(  H\right)
>0\ $is decreasing on some small interval $[0,\varepsilon),$ $\varepsilon>0.$
However, by $\lim_{H\downarrow0}w\left(  H\right)  =0$ we get a contradiction
and must have $w\left(  H\right)  \equiv0.$ Thus the two solutions are in fact equal.

\begin{remark}
Since the trajectory $\left(  H\left(  \xi\right)  ,G\left(  \xi\right)
\right)  \ $in $W^{u}\left(  O\right)  $ is unique, in the above
existence\ proof we can choose $\lambda_{1}$ as close to $p/c$ as possible. In
particular $\left(  H\left(  \xi\right)  ,G\left(  \xi\right)  \right)  $ in
$W^{u}\left(  O\right)  \ $must satisfy the following
\begin{equation}
\frac{p}{c}=\lim_{\xi\rightarrow-\infty}\frac{G\left(  \xi\right)  }{H\left(
\xi\right)  }=\lim_{\xi\rightarrow-\infty}\frac{pU^{p-1}\left(  \xi\right)
U^{\prime}\left(  \xi\right)  }{U^{p}\left(  \xi\right)  }=\lim_{\xi
\rightarrow-\infty}\frac{pU^{\prime}\left(  \xi\right)  }{U\left(  \xi\right)
}, \label{p/c}%
\end{equation}
which implies the asymptotic behavior%
\begin{equation}
\lim_{\xi\rightarrow-\infty}\frac{U^{\prime}\left(  \xi\right)  }{U\left(
\xi\right)  }=\frac{1}{c}. \label{1/c}%
\end{equation}

\end{remark}

\ \ \ \ \ \ \ 

Let $\left(  H_{c}\left(  \xi\right)  ,G_{c}\left(  \xi\right)  \right)  $
denote the trajectory whose existence and uniqueness have been established and
let $U_{c}\left(  \xi\right)  $ be the corresponding function of $\xi.$ That
is, $U_{c}\left(  \xi\right)  =H_{c}^{1/p}\left(  \xi\right)  .\ $Recall that
along any positive\ solution $U\left(  \xi\right)  $ of (\ref{Q1}), we
have$\ dE/d\xi=-2c\left(  U^{\prime}\right)  ^{2}/U^{p}\leq0,\ $%
where$\ E\left(  \xi\right)  $ is given by (\ref{E}). It follows that the
quantity$\ $%
\[
E\left(  H,G\right)  =\frac{G^{2}}{p^{2}}H^{2/p-2}+H^{2/p}-\frac{2}%
{2-p}H^{2/p-1},\ \ \ H=U^{p},\ \ \ G=H_{\xi}%
\]
is strictly decreasing on orbits$\ $of $X_{c},\ $except$\ c=0$\ (when\ $c=0$%
,\ all orbits are closed\ curves).$\ $Thus $H^{2/p}-2\left(  2-p\right)
^{-1}H^{2/p-1},$ and therefore $H,$ are bounded from above on any orbit of
$X_{c}.$ Furthermore, it also implies that $\left\vert G\right\vert $ is bounded.

Using the fact that $\left(  1,0\right)  $ is an attracting spiral point, one
can show that any orbit (here we only care about those orbits with positive
$H\ $everywhere)\ converges to $\left(  1,0\right)  ,$ and winds around this
point infinitely many times. In particular, any orbit will intersect the
$H$-axis (to see this, just look at the vector field (\ref{X})). For the
function $U_{c}\left(  \xi\right)  ,$ this means that it will converge to $1$
as $\xi\rightarrow\infty$ and that it will oscillate infinitely often around
its limit value. Its derivative $U_{c}^{\prime}\left(  \xi\right)  $ must
therefore vanish infinitely often; by replacing $U_{c}\left(  \xi\right)  $ by
$U_{c}\left(  \xi-\xi_{0}\right)  $ for some $\xi_{0}\in\mathbb{R}$ if
necessary, we may assume that the first zero of $U_{c}^{\prime}\ $is $\xi=0$
and $U_{c}^{\prime}\left(  \xi\right)  >0$ for all $\xi\in(-\infty,0).$

So far we have constructed the solution $U_{c}\left(  \xi\right)  $ satisfying
(i) and (ii) of (\ref{Q2}). Since we also have (\ref{1/c}), one can choose a
large constant$\lambda_{c}>0\ $so that$\,$(iii) of (\ref{Q2}) is also
satisfied for such $U_{c}\left(  \xi\right)  .$\newline

\ \ \ \ \ 

This complete our construction of $U_{c}.\ $To finishes the proof, we need to
verify (\ref{Q3}) for $U_{c}\left(  \xi\right)  .$

$\ \ \ \ \ \ \ $

We observe that the segment of $W^{u}\left(  O\right)  $ which lies in the
first quadrant is the graph of some function $G=g_{c}\left(  H\right)  $ for
$0\leq H\leq h_{c},$ where $\left(  h_{c},0\right)  $ is the first point of
intersection of $W^{u}\left(  O\right)  $ with the $H$-axis. Since
$U_{c}\left(  0\right)  =h_{c}^{1/p},$ we have to show that $h_{c}$\textbf{
}is\emph{ monotone decreasing} in\textbf{ }$c.$\textbf{ }

Let $c^{\prime}<c\ $be given, and suppose that $h_{c^{\prime}}\leq h_{c}.\ $We
want to derive a contradiction. Assume first that $h_{c^{\prime}}<h_{c}%
.\ $Compare the two vector fields $X_{c}\ $and $X_{c^{\prime}}\ $in the first
quadrant. If the backward orbit of $X_{c}\ $through $\left(  h_{c},0\right)
\ $and the backward orbit of $X_{c^{\prime}}\ $through $\left(  h_{c^{\prime}%
},0\right)  \ $intersect at a first point $\left(  H_{0},G_{0}\right)
,$$G_{0}>0,\ $we have the following comparison of the two vector fields at
$\left(  H_{0},G_{0}\right)  :$%
\[
-pH_{0}^{2}+pH_{0}-cG_{0}+\frac{p-1}{p}G_{0}^{2}<-pH_{0}^{2}+pH_{0}-c^{\prime
}G_{0}+\frac{p-1}{p}G_{0}^{2}%
\]
which implies that the backward orbit of $X_{c^{\prime}}\ $through $\left(
h_{c^{\prime}},0\right)  \ $cannot pass through the graph of $g_{c}\left(
H\right)  .\ $As a consequence, the graph of $g_{c^{\prime}}\left(  H\right)
\ $on the domain $0<H<h_{c^{\prime}}\ $must be below the graph of
$g_{c}\left(  H\right)  .\ $On the other hand by (\ref{1/c})\ we know%
\begin{equation}
g_{c}^{\prime}\left(  0\right)  =\lim_{\xi\rightarrow-\infty}\frac
{g_{c}\left(  H\left(  \xi\right)  \right)  }{H\left(  \xi\right)  }=\lim
_{\xi\rightarrow-\infty}\frac{G\left(  \xi\right)  }{H\left(  \xi\right)
}=\frac{p}{c} \label{gc1}%
\end{equation}
and similarly $g_{c^{\prime}}^{\prime}\left(  0\right)  =p/c^{\prime}.\ $Hence
$g_{c^{\prime}}^{\prime}\left(  0\right)  >g_{c}^{\prime}\left(  0\right)
\ $and this gives a contradiction%
\begin{equation}
\lim_{H\rightarrow0^{+}}g_{c^{\prime}}^{\prime}\left(  H\right)  =\frac
{p}{c^{\prime}}>\lim_{H\rightarrow0^{+}}g_{c}^{\prime}\left(  H\right)
=\frac{p}{c}>0. \label{g}%
\end{equation}
Thus $h_{c^{\prime}}<h_{c}\ $is impossible.\ \ 

If $c^{\prime}<c\ $but $h_{c^{\prime}}=h_{c},\ $then by continuity we must
have$\ 0<g_{c^{\prime}}\left(  H\right)  \leq g_{c}\left(  H\right)  \ $for
all\ $0<H<h_{c^{\prime}}.\ $But now estimate (\ref{g}) still holds and we
obtain the same contradiction.

\begin{remark}
By (\ref{gc1})\ and (\ref{1/c}), we have%
\begin{align*}
\frac{p}{c}  &  =g_{c}^{\prime}\left(  0\right)  =\lim_{H\rightarrow0^{+}%
}g_{c}^{\prime}\left(  H\right)  =\lim_{\xi\rightarrow-\infty}\frac{G_{\xi}%
}{H_{\xi}}\\
&  =\lim_{\xi\rightarrow-\infty}\frac{pU^{p-1}\left(  \xi\right)
U^{\prime\prime}\left(  \xi\right)  +p\left(  p-1\right)  U^{p-2}\left(
\xi\right)  \left(  U^{\prime}\left(  \xi\right)  \right)  ^{2}}%
{pU^{p-1}\left(  \xi\right)  U^{\prime}\left(  \xi\right)  }=\lim
_{\xi\rightarrow-\infty}\frac{U^{\prime\prime}\left(  \xi\right)  }{U^{\prime
}\left(  \xi\right)  }+\frac{p-1}{c}%
\end{align*}
and derive the limit
\begin{equation}
\lim_{\xi\rightarrow-\infty}\frac{U^{\prime\prime}\left(  \xi\right)
}{U^{\prime}\left(  \xi\right)  }=\frac{1}{c}. \label{2/c}%
\end{equation}
From\ (\ref{1/c})\ and (\ref{2/c}),\ it is not hard to see that asymptotically
$U\left(  \xi\right)  $ is given by $ae^{\left(  1/c\right)  \xi}\ $as
$\xi\rightarrow-\infty\ $for some constant $a>0.\ $
\end{remark}

\ \ \ \ 

A similar argument also show that $g_{c}\left(  H\right)  $ is a strictly
decreasing function of $c$ for fixed $H,$ i.e., as $c\downarrow0,$ the
unstable set $W_{c}^{u}\left(  O\right)  $ moves upwards.

We next claim that $U_{c}\left(  0\right)  \rightarrow\infty\ $%
as\ $c\downarrow0.\ $Assume that $U_{c}\left(  0\right)  $ were bounded, as
$c\downarrow0.$ Then the $h_{c}^{\prime}s$ would converge to some $h_{0}>1.$
The vector field $X_{c}$ is well-defined and smooth for all $c\in\mathbb{R},$
so the unstable set $W_{c}^{u}\left(  O\right)  ,$ being the orbit of $X_{c}$
through $\left(  h_{c},0\right)  ,$ must converge to the orbit of $X_{0}$
through $\left(  h_{0},0\right)  ,$ where
\[
X_{0}\left(  H,G\right)  =\left(  HG,\ -pH^{2}+pH+\frac{p-1}{p}G^{2}\right)
\]
and the quantity
\[
E\left(  H,G\right)  =\frac{G^{2}}{p^{2}}H^{2/p-2}+H^{2/p}-\frac{2}%
{2-p}H^{2/p-1},\ \ \ (H=U^{p},\ G=pU^{p-1}U_{\xi})
\]
is constant on the orbits of $X_{0}.$ By (\ref{ODE})\ we know that all orbits
of $X_{0}$ are periodic (\emph{due to the condition }$p>2$). In particular,
for $c=0,$ the orbit of $X_{0}$ through $\left(  h_{0},0\right)  $ will
satisfy the equation%
\[
\frac{G^{2}}{p^{2}}H^{2/p-2}+H^{2/p}-\frac{2}{2-p}H^{2/p-1}=h_{0}^{2/p}%
-\frac{2}{2-p}h_{0}^{2/p-1}>0,\ \ \ \text{ \ }h_{0}>1
\]
and from this equation we see that the orbit of $X_{0}$ through $\left(
h_{0},0\right)  $ will intersect the $H$-axis at some point $\left(  h_{\ast
},0\right)  ,$ $0<h_{\ast}<1,$ when followed backwards in time, where
$h_{\ast}\ $satisfies%
\[
h_{\ast}^{2/p}-\frac{2}{2-p}h_{\ast}^{2/p-1}=h_{0}^{2/p}-\frac{2}{2-p}%
h_{0}^{2/p-1}.
\]

By continuous dependence on parameters, the same will be true for some small
$c>0,$ a contradiction. Therefore we have $\lim_{c\downarrow0}U_{c}\left(
0\right)  =\infty.$

Recall that for fixed $c>0,$ the quantity
\[
E\left(  H\left(  \xi\right)  ,G\left(  \xi\right)  \right)  =\frac
{G^{2}\left(  \xi\right)  }{p^{2}}H^{2/p-2}\left(  \xi\right)  +H^{2/p}\left(
\xi\right)  -\frac{2}{2-p}H^{2/p-1}\left(  \xi\right)
\]
is strictly decreasing along the unstable orbit. We already know that as
$c\downarrow0,$ $h_{c}\left(  0\right)  \uparrow\infty.$ In
particular,\ $E\left(  H\left(  \xi\right)  ,G\left(  \xi\right)  \right)  $
is uniformly large on $(-\infty,0]$ since
\[
E\left(  H\left(  0\right)  ,G\left(  0\right)  \right)  =h_{c}^{2/p}\left(
0\right)  -\frac{2}{2-p}h_{c}{}^{2/p-1}\left(  0\right)  \rightarrow
\infty\ \ \ \text{as \ \ }c\downarrow0.
\]
Now when we confine to the region $\delta\leq H=U^{p}\leq\delta^{-1},$ as
$c\downarrow0,$ we must have $g_{c}\left(  H\right)  =G\uparrow\infty$ as
$c\downarrow0,$ and uniformly so on the interval $\delta\leq H\leq1/\delta.$\ 

Since$\ g_{c}\left(  H\right)  =G=pU_{c}^{p-1}U_{c}^{\prime}$ and $\delta\leq
H\leq1/\delta,$ we must have $U_{c}^{\prime}$ sufficiently large as
$c\downarrow0.$ Therefore (\ref{Q3}) also holds.\ The proof of Theorem
\ref{thm-type-one-1}\ is done.$%
\hfill
\square\ $

\subsection{Proof of Theorem \ref{thm-type-one}.}

Assume $v\left(  x,t\right)  $ is a \emph{type-one solution} to equation
$\left(  \clubsuit\right)  $\ with $p\geq2$. Then the rescaled
positive\ function $u\left(  x,\tau\right)  $ is bounded from above.\ Choose a
large constant $A\ $so that%
\begin{equation}
\left\{
\begin{array}
[c]{l}%
u\left(  x,\tau\right)  \leq A\ \ \ \text{for all\ \ \ }\left(  x,\tau\right)
\in S_{m}^{1}\times\lbrack0,\infty)%
\vspace{3mm}%
\\
\left\vert u_{x}\left(  x,0\right)  \right\vert \leq A\ \ \ \text{for
all\ \ \ }x\in S_{m}^{1}%
\vspace{3mm}%
\\
u\left(  x,0\right)  \geq\frac{1}{A}\ \ \ \text{for all\ \ \ }x\in S_{m}^{1}.
\end{array}
\right.  \label{Q4}%
\end{equation}
Also choose $c>0$ so small that the solution $U_{c}\left(  \xi\right)  $ of
the last section satisfies$\ U_{c}\left(  0\right)  >A\ $and$\ U_{c}^{\prime
}\left(  \xi\right)  >A$ whenever$\ A^{-1}\leq U_{c}\left(  \xi\right)  \leq
A.\ $By Theorem\ \ref{thm-type-one-1} such\ a $c$ exists, together with the
existence of a large constant $\lambda_{c}>0$ such that $0<U_{c}^{\prime
}\left(  \xi\right)  \leq\lambda_{c}U_{c}\left(  \xi\right)  \ $for
all\ $\xi\in(-\infty,0].\ $Note that here the number $c\ $and $\lambda_{c}$
both depend on the initial data $u\left(  x,0\right)  .\ $

For any$\ $fixed$\ \left(  x_{0},\tau_{0}\right)  $ we have $0<u\left(
x_{0},\tau_{0}\right)  <U_{c}\left(  0\right)  $ and since $U_{c}\left(
\xi\right)  $ is strictly increasing on$\ (-\infty,0]$ there exists a unique
$x_{1}\in\mathbb{R}$ for which$\ U_{c}\left(  x_{1}-c\tau_{0}\right)
=u\left(  x_{0},\tau_{0}\right)  .$ Consider the function
\[
u^{\ast}\left(  x,\tau\right)  =U_{c}\left(  x-x_{0}+x_{1}-c\tau\right)
,\ \ \ u^{\ast}\left(  x,\tau\right)
\]
then$\ u^{\ast}\left(  x,\tau\right)  $ is a solution of$\ $the
equation$\ \partial u/\partial\tau=u^{p}\left(  u_{xx}+u-u^{1-p}\right)  $ on
the region%
\[
Q=\left\{  \left(  x,\tau\right)  :x<x_{0}-x_{1}+c\tau,\ \tau>0\right\}
\]
and the difference
\begin{equation}
w\left(  x,\tau\right)  =u^{\ast}\left(  x,\tau\right)  -u\left(
x,\tau\right)  =U_{c}\left(  x-x_{0}+x_{1}-c\tau\right)  -u\left(
x,\tau\right)
\end{equation}
satisfies a linear parabolic PDE of the form (see \cite{ANG}, p. 608)%
\begin{equation}
\frac{\partial w}{\partial\tau}=a\left(  x,\tau\right)  w_{xx}+b\left(
x,\tau\right)  w_{x}+c\left(  x,\tau\right)  w. \label{linear}%
\end{equation}
We note that $w\left(  x_{0},\tau_{0}\right)  =0$ and on the
boundary$\ \partial Q\bigcap\left\{  \tau>0\right\}  $ (i.e.,\ when
$x=x_{0}-x_{1}+c\tau$)\ we have%
\[
w\left(  x,\tau\right)  =U_{c}\left(  0\right)  -u\left(  x,\tau\right)  \geq
U_{c}\left(  0\right)  -A>0.
\]
On the other part of $\partial Q,\ $i.e.,\ when $x<x_{0}-x_{1}\ $and $\tau=0$
we have%
\[
w\left(  x_{0}-x_{1},0\right)  =U_{c}\left(  0\right)  -u\left(  x_{0}%
-x_{1},0\right)  >0
\]
and$\ w\left(  x,0\right)  =U_{c}\left(  x-x_{0}+x_{1}\right)  -u\left(
x,0\right)  $ becomes negative as $x\rightarrow-\infty$ due to (\ref{Q4}).
Hence $w\left(  x,0\right)  $ must have at least one zero $y_{0}\ $on $\left(
-\infty,x_{0}-x_{1}\right)  .\ $At any$\ $zero$\ y_{0}$ we have$\ $%
\[
\frac{1}{A}\leq u\left(  y_{0},0\right)  =U_{c}\left(  y_{0}-x_{0}%
+x_{1}\right)  \leq A
\]
and so
\[
w_{x}\left(  y_{0},0\right)  =U_{c}^{\prime}\left(  y_{0}-x_{0}+x_{1}\right)
-u_{x}\left(  y_{0},0\right)  >A-u_{x}\left(  y_{0},0\right)  \geq0.
\]
Hence$\ w\left(  x,0\right)  $ cannot have more than one zero on the
interval$\ \left(  -\infty,x_{0}-x_{1}\right)  .\ $

By the Sturmian theorem, the number of zeros of $x\rightarrow w\left(
x,\tau\right)  ,$ counted with multiplicity,\ cannot increase with time.\ Now
by our construction we have$\ w\left(  x_{0},\tau_{0}\right)  =0\ $and
since\ this is \emph{the only zero} of$\ w\left(  \cdot,\tau_{0}\right)
,\ $we must have$\ w_{x}\left(  x_{0},\tau_{0}\right)  >0\ $(since $w\left(
x_{0}-x_{1}+c\tau_{0},\tau_{0}\right)  >0$ and $w\left(  -\infty,\tau
_{0}\right)  <0$).\ Thus%
\begin{equation}
u_{x}\left(  x_{0},\tau_{0}\right)  <U_{c}^{\prime}\left(  x_{1}-c\tau
_{0}\right)  \leq\lambda_{c}U_{c}\left(  x_{1}-c\tau_{0}\right)  =\lambda
_{c}u\left(  x_{0},\tau_{0}\right)  .
\end{equation}

By applying the same argument to$\ u\left(  -x,\tau\right)  $ one can also
obtain $-u_{x}\leq\lambda_{c}u,$ so that $\left\vert u_{x}\left(
x,\tau\right)  \right\vert \leq\lambda_{c}u\left(  x,\tau\right)  $ for all
$\left(  x,\tau\right)  \in\in S_{m}^{1}\times\lbrack0,\infty).\ $The proof of
Theorem \ref{thm-type-one}\ is done.$%
\hfill
\square$

\subsection{Proof of type-one convergence.}

To go further we need to look more closely at the following\ ODE:%
\begin{equation}
w^{\prime\prime}\left(  x\right)  +w\left(  x\right)  -w^{1-p}\left(
x\right)  =0,\text{\ \ }\ x\in\left(  -\infty,\infty\right)  ,\ \ \ p>2.
\label{W-ODE}%
\end{equation}
It is easy to see that any solution $w\left(  x\right)  $ to it is
\emph{positive everywhere} and\emph{ periodic} over $x\in\left(
-\infty,\infty\right)  \ $(this property is valid for $p\geq2;\ $when
$p\in\left(  0,2\right)  ,\ w\left(  x\right)  $ may have different behavior,
see \cite{LPT}\ and \cite{PT})$.\ $Let $a\leq1$ be the minimal value
of$\ w\left(  x\right)  $ on$\ \left(  -\infty,\infty\right)  $.\ Without loss
of generality, we may assume that $a=w\left(  0\right)  \ $(and so $w^{\prime
}\left(  0\right)  =0$)$\ $and by reflection$\ $(if $w\left(  x\right)  $ is a
solution, so is $w\left(  -x\right)  $)$\ w\left(  x\right)  $ must be
symmetric with respect to any local maximum point or minimum point. It also
satisfies the energy identity
\begin{equation}
\left(  w^{\prime}\left(  x\right)  \right)  ^{2}+w^{2}\left(  x\right)
-\frac{2}{2-p}w^{2-p}\left(  x\right)  =F\left(  a\right)
\ \ \ \text{for\ all}\ \ \ x\in\left(  -\infty,\infty\right)  \label{wF}%
\end{equation}
where$\ F\left(  a\right)  =a^{2}-2\left(  2-p\right)  ^{-1}a^{2-p}>0.\ $For
$p>2,\ $the$\ $convex\ positive\ function $F\left(  s\right)  =s^{2}-2\left(
2-p\right)  ^{-1}s^{2-p}\ $decreases on $s\in\left(  0,1\right)
\ $with$\ \lim_{s\rightarrow0^{+}}F\left(  s\right)  =+\infty,$\ and increases
to $+\infty\ $on $\left(  1,\infty\right)  $. Given $a\in(0,1],\ $there is a
unique $b\geq1\ $so that $F\left(  a\right)  =F\left(  b\right)  ,\ $where
$b=\max_{x\in\mathbb{R}}w\left(  x\right)  ,\ $and the minimal period
$T=2R\left(  a\right)  \ $of $w\left(  x\right)  $ is given by%
\begin{equation}
T=2R\left(  a\right)  =2\int_{a}^{b}\frac{ds}{\sqrt{F\left(  a\right)
-F\left(  s\right)  }}=2\int_{a}^{b}\frac{ds}{\sqrt{\left(  a^{2}-\frac
{2}{2-p}a^{2-p}\right)  -\left(  s^{2}-\frac{2}{2-p}s^{2-p}\right)  }%
},\ \ \ F\left(  b\right)  =F\left(  a\right)  . \label{T}%
\end{equation}
The above integral is improper near both$\ a$ and $b.\ $

It has been shown in Urbas \cite{U2} that
\begin{equation}
\lim_{a\rightarrow0^{+}}R\left(  a\right)  =\frac{\pi}{2},\ \ \ \lim
_{a\rightarrow1^{-}}R\left(  a\right)  =\frac{\pi}{\sqrt{p}},\ \ \ p\in\left(
2,\infty\right)  . \label{ur}%
\end{equation}
Moreover, by Corollary 5.6 of Andrews\ \cite{AN3},\ we know that $R\left(
a\right)  $ is strictly decreasing in $a\in\left(  0,1\right)  $ when
$p\in\left(  4,\infty\right)  $\ and strictly increasing in $a\in\left(
0,1\right)  $ when $p\in\left(  2,4\right)  .\ $When $p=4,$ all solutions of
equation (\ref{W-ODE}) are $\pi$-periodic (see (\ref{wb}) also).\ 

\begin{remark}
\label{rmk1}As a comparison, when $p\in\left(  0,2\right)  $ we have\ (see
\cite{U2} and \cite{AN3} again)%
\[
\lim_{a\rightarrow0^{+}}R\left(  a\right)  =\frac{\pi}{p},\ \ \ \lim
_{a\rightarrow1^{-}}R\left(  a\right)  =\frac{\pi}{\sqrt{p}},\ \ \ p\in\left(
0,2\right)
\]
and $R\left(  a\right)  $ is strictly decreasing in $a\in\left(  0,1\right)  $
when $p\in\left(  0,1\right)  $ and strictly increasing in $a\in\left(
0,1\right)  $ when $p\in\left(  1,2\right)  .\ $When $p=1,\ $all solutions to
the ODE (\ref{W-ODE}) has period $2\pi.\ $
\end{remark}

One can also write the ODE (\ref{W-ODE})\ as a system
\begin{equation}
\dfrac{dw}{dx}=h,\ \ \ \dfrac{dh}{dx}=-w+w^{1-p},\ \ \ p>2.
\end{equation}
Then the vector field $V\left(  w,h\right)  =\left(  h,-w+w^{1-p}\right)  $
has only one equilibrium point $\left(  1,0\right)  $ on the
half-plane$\ \left\{  w>0\right\}  \ $and the eigenvalues of the linearization
at it are $\lambda=\pm\sqrt{p}i\ $(this matches with the second limit of
(\ref{ur})).\ The phase portrait of $V\ $on $\left\{  w>0\right\}  \ $is a
family of \emph{closed orbits} $C\left(  a\right)  \ $centered at$\ \left(
1,0\right)  \ $with period $2R\left(  a\right)  ,\ $where $a=\min
_{x\in\mathbb{R}}w\left(  x\right)  $.\ Thus the intersections of $C\left(
a\right)  $ and the $w$-axis are $\left(  a,0\right)  \ $and $\left(
b,0\right)  \ $with$\ F\left(  a\right)  =F\left(  b\right)  ,$ $w\left(
0\right)  =a\leq1,\ w\left(  R\left(  a\right)  \right)  =b\geq1,\ w^{\prime
}\left(  0\right)  =w^{\prime}\left(  R\left(  a\right)  \right)  =0.$

We can now state the following convergence theorem, which is a generalization
of Theorem A of \cite{ANG}:

\begin{theorem}
\label{thmA}(\emph{convergence of type-one blow-up for}\textsf{ }$p>2$)\ Let
$p>2$ and let $v\left(  x,t\right)  >0\ $be a \emph{type-one solution} of
$\left(  \clubsuit\right)  \ $defined on some maximal time interval
$[0,T_{\max})$.\ Then as $\tau\rightarrow\infty$ the rescaled solution
$u\left(  x,\tau\right)  ,$ given by (\ref{rescale}),\ converges in
$C^{\infty}\left(  S_{m}^{1}\right)  $ to a smooth\ positive $2m\pi$-periodic
function $w\left(  x\right)  ,\ $which is an entire solution of the ODE
\begin{equation}
w^{\prime\prime}\left(  x\right)  +w\left(  x\right)  -w^{1-p}\left(
x\right)  =0\ \ \ \text{for all\ \ }\ x\in\mathbb{R}. \label{wode}%
\end{equation}

\end{theorem}%

\proof
Assume type-one blow-up of $v\left(  x,t\right)  $.\ For any sequence
$\tau_{n}\rightarrow\infty$ by Arzela-Ascoli theorem there is a subsequence,
which\ we also call it $\tau_{n}$, so that $u\left(  x,\tau_{n}\right)  $
converges uniformly on $S_{m}^{1}$ to a Lipschitz\ function $w\left(
x\right)  \geq0$, which is $2m\pi$-periodic. By\ Theorem \ref{thm-type-one}%
$,\ u\left(  x,\tau\right)  \ $has positive lower bound $e^{-2\lambda m\pi}$
for all $\tau\ $(see (\ref{KK1})),$\ $hence $w\left(  x\right)  \ $is strictly
positive everywhere.\ Now we can apply similar argument as in Proposition 12
and 14 of \cite{LPT} (since $p>2,$\ $u\left(  x,\tau\right)  \ $has positive
lower bound is essential in (32),\ p.160 of \cite{LPT}) to obtain the
conclusion that$\ w\left(  x\right)  $ satisfies the ODE (\ref{wode})
everywhere.\ By regularity theory for uniform parabolic equations,$\ w\left(
x\right)  \ $is smooth\ and we have $C^{\infty}$ convergence of $u\left(
x,\tau_{n}\right)  $ to $w\left(  x\right)  \ $as $\tau_{n}\rightarrow\infty$.

If we does not have full time convergence of $u\left(  x,\tau\right)  $
as$\ \tau\rightarrow\infty,\ $then there will exist two sequence of times
$\tau_{n}\rightarrow\infty$ and$\ \tilde{\tau}_{n}\rightarrow\infty$ such that
$u\left(  x,\tau_{n}\right)  \rightarrow w\left(  x\right)  $ and $u\left(
x,\tilde{\tau}_{n}\right)  \rightarrow\tilde{w}\left(  x\right)  ,\ $where
$w,$ $\tilde{w}\ $are \emph{different}\ positive $2m\pi$-periodic solutions of
the ODE (\ref{wode}).\ Let $a,\ \tilde{a}\in(0,1]$ be the minimum values of
$w,\ \tilde{w}.\ $We may assume $a\leq\tilde{a}.\ $Note that although
$w\left(  x\right)  \ $is different from $\tilde{w}\left(  x\right)  $,
it\ may be possible to have $a=\tilde{a}$.\ By\ the above discussion, we
have$\ w\left(  R\left(  a\right)  \right)  =b\geq w\left(  R\left(  \tilde
{a}\right)  \right)  =\tilde{b}\in\lbrack1,\infty).$

If we have $a=\tilde{a},$ then $w\left(  x\right)  $ must be a translation of
$\tilde{w}\left(  x\right)  $ and we can find some $x_{0}\in\mathbb{R}$ with
$w^{\prime}\left(  x_{0}\right)  \tilde{w}^{\prime}\left(  x_{0}\right)
<0\ $(i.e., they have different\ signs). This would contradict Proposition 23
of \cite{LPT}. Therefore we only have to consider the case $a<\tilde{a}.$

\ \ \ \ \ 

For $a<\tilde{a}$ there are two cases to discuss.\ 

\ \ 

\underline{Case 1}:\ $p\in\left(  2,\infty\right)  ,\ p\neq4.$

\ \ \ \ \ 

By the discussion before Remark \ref{rmk1}, in such case we must have
$R\left(  a\right)  \neq R\left(  \tilde{a}\right)  \ $since $R\left(
a\right)  $ is a monotone function in $a\in\left(  0,1\right)  .\ $Now
$w\left(  x\right)  $ and $\tilde{w}\left(  x\right)  $ have
different\ periods and we can find some $x_{0}\in\mathbb{R}$ such that
$w^{\prime}\left(  x_{0}\right)  \tilde{w}^{\prime}\left(  x_{0}\right)
<0.\ $We obtain the same contradiction due to Proposition 23 of \cite{LPT}.\ 

\ \ \ \ 

\underline{Case 2}:\ $p=4.$

\ \ \ \ 

In this case by (\ref{wb}),\ up to a translation, all solutions to the ODE
$w^{\prime\prime}+w-w^{-3}=0$ are $\pi$-periodic (see \cite{AN3}%
,\ \cite{U2})\ and are given by (if $w(0)\geq1\ $is the maximum) the
$1$-parameter family of functions in (\ref{wb}). Unfortunately now for any
$x_{0}\in\mathbb{R}$ we have$\ w^{\prime}\left(  x_{0}\right)  \tilde
{w}^{\prime}\left(  x_{0}\right)  \geq0\ $and thus Proposition 23 of
\cite{LPT} is not applicable here.\ A different\ method has to be used
here\footnote{D.-H. Tsai would like to thank Prof. Matano\ for teaching him
the zero-number argument several years ago.\ It is now used in the proof of
Theorem \ref{thmA}.}.

$\ $Let $\mathcal{Z}\left[  w-\tilde{w}\right]  $ denote the number of
zero\ points $\xi\in S_{m}^{1}\ $(or $\xi\in\lbrack-m\pi,m\pi)$)
with$\ w\left(  \xi\right)  -\tilde{w}\left(  \xi\right)  =0.$ Also let $Y$
denote the function space of all solutions of the ODE (\ref{wode}) on
$\mathbb{R\ }$(since now $p=4,\ $all solutions to the ODE (\ref{wode}) has
minimal period $\pi$)$.\ $For any $\xi\in S_{m}^{1}$ with $w\left(
\xi\right)  -\tilde{w}\left(  \xi\right)  =0$, by uniqueness\ we must have
$w^{\prime}\left(  \xi\right)  \neq\tilde{w}^{\prime}\left(  \xi\right)  .$
Hence at each intersection point the graphs of the two functions
$w,\ \tilde{w}\ $are transversal.\ 

For any $z\left(  x\right)  \in Y,\ $the difference $u\left(  x,\tau\right)
-z\left(  x\right)  $ satisfies a linear parabolic equation of the form
(\ref{linear}) (since $z\left(  x\right)  $ is also a solution to the PDE
(\ref{dudt})). By Angenent's result in p. 607 of \cite{ANG}\ (Lemma 2.4\ in
p.\ 165 of Chen-Matano \cite{CM} is more applicable here), the number
$\mathcal{Z}\left[  u\left(  \cdot,\tau\right)  -z\left(  \cdot\right)
\right]  $ is \emph{non-increasing} in time $\tau\in\left(  0,\infty\right)
.\ $Also note that we have the convergence of $u\left(  x,\tau_{n}\right)
\ $to $w\left(  x\right)  $ in $C^{1},$ which implies
\[
\mathcal{Z}\left[  u\left(  \cdot,\tau_{n}\right)  -z\left(  \cdot\right)
\right]  =\mathcal{Z}\left[  w-\tilde{w}\right]
\]
for all large $n$ and all $z\in Y$ that are sufficiently close to $\tilde{w}$
in $C^{1}$ norm on $S_{m}^{1}.\ $In particular, we can conclude the
following:\ there exists a time $T>0\ $and a number $\delta>0\ $such that
\begin{equation}
\mathcal{Z}\left[  u\left(  \cdot,\tau\right)  -z\left(  \cdot\right)
\right]  =\mathcal{Z}\left[  w-\tilde{w}\right]  \label{co}%
\end{equation}
for all $\tau>T$ and all $z\in Y$ satisfying $\left\Vert z-\tilde
{w}\right\Vert _{C^{1}\left(  S_{m}^{1}\right)  }<\delta.\ $

The number $\mathcal{Z}\left[  u\left(  \cdot,\tau\right)  -z\left(
\cdot\right)  \right]  $ remains a constant for large time implies that the
function $x\rightarrow u\left(  x,\tau\right)  -z\left(  x\right)  $ does not
have a degenerate zero (i.e., multiple zero)\ in $S_{m}^{1}$ for any fixed
$\tau>T\ $(see \cite{CM}).$\ $But since $u\left(  x,\tilde{\tau}_{n}\right)  $
converges to $\tilde{w}\left(  x\right)  \ $in\ $C^{1}\left(  S_{m}%
^{1}\right)  $ norm as $n\rightarrow\infty,$ the graph of the function
$x\rightarrow u\left(  x,\tilde{\tau}_{n}\right)  $ must be tangential to the
graph of some $z\in Y\ $satisfying $\left\Vert z-\tilde{w}\right\Vert
_{C^{1}\left(  S_{m}^{1}\right)  }<\delta.\ $For example, for fixed $x_{0}%
\ $one can choose $z\left(  x\right)  $ to be the solution of
\begin{equation}
\left\{
\begin{array}
[c]{l}%
z^{\prime\prime}\left(  x\right)  +z\left(  x\right)  -z^{-3}\left(  x\right)
=0%
\vspace{3mm}%
\\
z\left(  x_{0}\right)  =u\left(  x_{0},\tilde{\tau}_{n}\right)
,\ \ \ z^{\prime}\left(  x_{0}\right)  =u_{x}\left(  x_{0},\tilde{\tau}%
_{n}\right)
\end{array}
\right.  \label{z}%
\end{equation}
then as $n$ large enough, $z\left(  x\right)  $ will be close to $\tilde
{w}\left(  x\right)  $ in $C^{1}\left(  S_{m}^{1}\right)  $ since $u\left(
x_{0},\tilde{\tau}_{n}\right)  \ $is close to $\tilde{w}\left(  x_{0}\right)
\ $and $u_{x}\left(  x_{0},\tilde{\tau}_{n}\right)  $ is close to $\tilde
{w}^{\prime}\left(  x_{0}\right)  .$ Now $u\left(  x,\tilde{\tau}_{n}\right)
-z\left(  x\right)  $ has a degenerate zero at $x_{0},\ $which is a
contradiction. \ 

\begin{remark}
Since $p=4,$\ $z\left(  x\right)  $ to the ODE (\ref{z}) has minimal period
$\pi.\ $In particular, it implies that $z\left(  x\right)  \in C^{1}\left(
S_{m}^{1}\right)  $.
\end{remark}

The above contradiction for either Case 1 or Case 2 implies that $w\left(
x\right)  \equiv\tilde{w}\left(  x\right)  $ and the proof is done.$%
\hfill
\square$

\ \ 

Theorem \ref{thmA}\ implies that for the contracting flow $\left(
\bigstar\right)  $, if $k_{\max}\left(  t\right)  \left(  T_{\max}-t\right)
^{1/\left(  \alpha+1\right)  }$ remains bounded as $t\rightarrow T_{\max},$
then the rescaled curvature%
\[
K\left(  x,\tau\right)  =\left(  p^{1/p}T_{\max}^{1/p}e^{-\tau}\right)
^{1/\alpha}k\left(  x,T_{\max}\left(  1-e^{-p\tau}\right)  \right)
,\ \ \ \alpha\in(0,1],\ \ \ \tau\in\lbrack0,\infty)
\]
converges in $C^{\infty}$ to a positive $K\left(  x\right)  \in C^{\infty
}\left(  S_{m}^{1}\right)  ,\ $which satisfies the ODE%
\begin{equation}
\left(  K^{\alpha}\right)  ^{\prime\prime}\left(  x\right)  +K^{\alpha}\left(
x\right)  -\frac{1}{K\left(  x\right)  }=0\ \ \ \text{for all\ \ \ }%
x\in\mathbb{R}. \label{kk-ode}%
\end{equation}
Geometrically this says that the evolving convex immersed closed curve
$\gamma_{t}$ shrinks to a point in an asymptotically self-similar way.\ 

\section{Type-two blow-up.}

We now turn to the much more difficult type-one blow-up. We point out that in
the proof of Theorem \ref{thm-type-one}, the integral condition
(\ref{integral-cond}) does not come into play at all. Hence even it is not
satisfied, Theorem \ref{thm-type-one} still holds.\ In view of this, we have
the following interesting observation:

\begin{lemma}
\label{lem3}(\textsf{existence of type-two blow-up for }$p\geq2$)\ Assume
$v_{0}\left(  x\right)  >0\in C^{\infty}\left(  S_{m}^{1}\right)  $ in
$\left(  \clubsuit\right)  $ does not satisfy (\ref{integral-cond}), i.e.,
\begin{equation}
\int_{S_{m}^{1}}v_{0}^{1-p}\left(  x\right)  e^{ix}dx\neq0 \label{not-zero}%
\end{equation}
then we have\ type-two blow-up\textsf{ }for the solution $v\left(  x,t\right)
$ to $\left(  \clubsuit\right)  $,\textsf{ }which means%
\begin{equation}
\limsup_{t\rightarrow T_{\max}}\left(  v_{\max}\left(  t\right)  \left(
T_{\max}-t\right)  ^{1/p}\right)  =\infty. \label{type-1}%
\end{equation}

\end{lemma}%

\proof
Without loss of generality we may assume
\[
\int_{S_{m}^{1}}v_{0}^{1-p}\left(  x\right)  \cos xdx>0.
\]
Since we have
\begin{equation}
\int_{S_{m}^{1}}v^{1-p}\left(  x,t\right)  e^{ix}dx=\int_{S_{m}^{1}}%
v_{0}^{1-p}\left(  x\right)  e^{ix}dx
\end{equation}
for all $t\in\lbrack0,T_{\max}),$ $u\left(  x,\tau\right)  $ satisfies$\ $%
\[
\lim_{\tau\rightarrow\infty}\int_{S_{m}^{1}}u^{1-p}\left(  x,\tau\right)  \cos
xdx=\lim_{\tau\rightarrow\infty}\left(  p^{1/p}T_{\max}^{1/p}e^{-\tau}\right)
^{1-p}\int_{S_{m}^{1}}v_{0}^{1-p}\left(  x\right)  \cos xdx=\infty
\]
which means that $\liminf_{\tau\rightarrow\infty}u_{\min}\left(  \tau\right)
=0.\ $If we have\ type-one blow-up, then Theorem \ref{thm-type-one} would
imply\ a positive lower bound of $u_{\min}\left(  \tau\right)  ,$ a
contradiction.$%
\hfill
\square$

\begin{remark}
Thus for $p\geq2,\ $\textsf{type-two blow-up} \emph{in equation}$\ \left(
\clubsuit\right)  $\emph{ is generic}.\ Moreover, type-one blow-up occurs only
when the initial data satisfies the integral condition (\ref{integral-cond}).
\end{remark}

\begin{remark}
When (\ref{integral-cond})\ is satisfied, then either type-one or type-two
blow-up can happen.\ For type-one, just take a separable solution of $\left(
\clubsuit\right)  $\ of the form $v\left(  x,t\right)  =h\left(  t\right)
g\left(  x\right)  ,$ where $g\left(  x\right)  >0$ on $S_{m}^{1}$ satisfies
the ODE\ (\ref{wode}) and $h\left(  t\right)  $ satisfies$\ dh/dt=h^{1+p}%
,\ h\left(  0\right)  >0.\ $For type-two, choose a convex immersed plane curve
with one big loop and one tiny loop. Then the corresponding evolution will
become singular without shrinking to a point in an asymptotically self-similar
way. Hence we obtain a type-two blow-up. The difficulty lies in the estimate
of blow-up rate.\ 
\end{remark}

\subsection{A special symmetric case for type-two blow-up and convergence.}

In this section we assume the initial data $v_{0}\left(  x\right)  >0$ to
equation $\left(  \clubsuit\right)  $ satisfies (\ref{integral-cond})\ and the
following \emph{symmetric condition }%
\begin{equation}
v_{0}\left(  x\right)  =v_{0}\left(  -x\right)  \ \ \ \text{and \ \ }%
v_{0}^{\prime}\left(  x\right)  <0,\ \ \ \ \ \text{for all\ \ \ }x\in\left(
0,m\pi\right)  . \label{v0}%
\end{equation}
If $v\left(  x,t\right)  $ is a solution to $\left(  \clubsuit\right)  $ with
the above\ initial data $v_{0}\left(  x\right)  $ then $\tilde{v}\left(
x,t\right)  :=v\left(  -x,t\right)  $ is also a solution to $\left(
\clubsuit\right)  $ with $\tilde{v}\left(  x,0\right)  =v_{0}\left(
-x\right)  =v_{0}\left(  x\right)  $ for all $x\in S_{m}^{1}.$ By uniqueness
we must have
\[
\ v\left(  x,t\right)  =v\left(  -x,t\right)  \ \ \ \text{for all\ \ \ }%
\left(  x,t\right)  \in\left(  0,m\pi\right)  \times\lbrack0,T_{\max})
\]
which also implies
\begin{equation}
v_{x}\left(  0,t\right)  =v_{x}\left(  m\pi,t\right)  =0\ \ \ \text{for
all\ \ \ }t\in\lbrack0,T_{\max}). \label{v1}%
\end{equation}
Also the second condition of (\ref{v0}) implies that $v_{0}^{\prime}\left(
x\right)  $ has exactly two zeros on $S_{m}^{1}$ and since the number of zeros
for $v_{x}\left(  x,t\right)  $ is nonincreasing in time, we must have%
\begin{equation}
v_{x}\left(  x,t\right)  <0\ \ \ \text{for all\ \ \ }\left(  x,t\right)
\in\left(  0,m\pi\right)  \times\lbrack0,T_{\max}). \label{decr}%
\end{equation}
Hence the two conditions of (\ref{v0}) are preserved for all time.\ In
particular, we have $v_{\max}\left(  t\right)  =v\left(  0,t\right)  \ $for
$t\in\lbrack0,T_{\max}).\ $

The main result in this section is the following convergence behavior for
type-two blow-up. One can view it as a partial generalization of Theorem C\ of
\cite{ANG} to the case $p\geq2$ since\ here we assume $v_{0}\left(  x\right)
\ $is symmetric\ and our convergence is only uniform,\ weaker than Angenent's
$C^{\infty}$ convergence. However, the advantage of focusing on the symmetric
case (\ref{v0})\ is that we always have \emph{type-two} blow-up and the proof
of convergence in Theorem \ref{thmC}\ below\ is very simple and straightforward.

\begin{lemma}
\label{lem4}Assume $v\left(  x,t\right)  \ $is a positive solution of $\left(
\clubsuit\right)  $ in $S_{m}^{1}$ (with $m\geq2$) where\ $v_{0}\left(
x\right)  $ satisfies (\ref{integral-cond}) and (\ref{v0}).\ Then$\ v\left(
x,t\right)  \ $has \textbf{type-two} blow-up.
\end{lemma}

%

\proof
Basically, we follow the arguments in p. 630 of \cite{ANG}.\ If $v\left(
x,t\right)  \ $has type-one blow-up, then by (\ref{KK1}) we have
\begin{equation}
v\left(  x,t\right)  \rightarrow\infty\ \text{as\ }t\rightarrow T_{\max
}\ \ \ \text{for all\ \ \ }x\in\left[  -m\pi,m\pi\right]  \text{.}
\label{v-go}%
\end{equation}
That is, the blow-up set of $v\left(  x,t\right)  \ $is the whole domain.\ 

When $m\ $is even,\ $m=2k,\ k\geq1,\ $consider the function%
\[
D\left(  t\right)  =\int_{0}^{\left(  2k-1\right)  \pi}v^{1-p}\left(
x,t\right)  \cos xdx.
\]
By (\ref{v-go}), we have
\begin{equation}
D\left(  t\right)  \rightarrow0\ \ \ \text{as}\ \ \ t\rightarrow T_{\max}.
\label{Dt}%
\end{equation}
Now by $\left(  \clubsuit\right)  $ we compute%
\begin{align*}
D^{\prime}\left(  t\right)   &  =\left(  1-p\right)  \int_{0}^{\left(
2k-1\right)  \pi}\left(  v_{xx}\left(  x,t\right)  +v\left(  x,t\right)
\right)  \cos xdx\\
&  =\left(  1-p\right)  \left(  v\left(  x,t\right)  \sin x+v_{x}\left(
x,t\right)  \cos x\right)  \mid_{0}^{\left(  2k-1\right)  \pi}=\left(
p-1\right)  v_{x}\left(  \left(  2k-1\right)  \pi,t\right)  <0
\end{align*}
due to (\ref{v1})\ and (\ref{decr}). Hence $D\left(  t\right)  $ is decreasing
and by (\ref{Dt}), it is positive for all $t\in\lbrack0,T_{\max}).\ $Also the
symmetry of $v\left(  x,t\right)  \ $implies that
\[
\int_{0}^{2k\pi}v^{1-p}\left(  x,t\right)  \cos xdx=\frac{1}{2}\int_{-2k\pi
}^{2k\pi}v^{1-p}\left(  x,t\right)  \cos xdx=0.
\]
Thus we have%
\[
0=\int_{0}^{\left(  2k-1\right)  \pi}v^{1-p}\left(  x,t\right)  \cos
xdx+\int_{\left(  2k-1\right)  \pi}^{2k\pi}v^{1-p}\left(  x,t\right)  \cos
xdx
\]
and then%
\[
\int_{\left(  2k-1\right)  \pi}^{2k\pi}v^{1-p}\left(  x,t\right)  \cos
xdx<0\ \ \ \text{for all\ \ \ }t\in\lbrack0,T_{\max}).
\]
However, by (\ref{decr})$\ $we have$\ $%
\[
\int_{\left(  2k-1\right)  \pi}^{2k\pi}v^{1-p}\left(  x,t\right)  \cos
xdx>0\ \ \ \text{for all\ \ \ }t\in\lbrack0,T_{\max}),
\]
which gives a contradiction.

When $m\ $is odd,\ $m=2k+1,\ k\geq1,\ $we consider the function%
\[
D\left(  t\right)  =\int_{0}^{2k\pi}v^{1-p}\left(  x,t\right)  \cos xdx
\]
and again by (\ref{v-go}), we have (\ref{Dt}). Now%
\begin{align*}
D^{\prime}\left(  t\right)   &  =\left(  1-p\right)  \int_{0}^{2k\pi}\left(
v_{xx}\left(  x,t\right)  +v\left(  x,t\right)  \right)  \cos xdx\\
&  =\left(  1-p\right)  \left(  v\left(  x,t\right)  \sin x+v_{x}\left(
x,t\right)  \cos x\right)  \mid_{0}^{2k\pi}=\left(  1-p\right)  v_{x}\left(
2k\pi,t\right)  >0
\end{align*}
and so $D\left(  t\right)  $ is increasing and therefore negative for all
time. By symmetry again we obtain%
\[
\int_{0}^{\left(  2k+1\right)  \pi}v^{1-p}\left(  x,t\right)  \cos
xdx=\frac{1}{2}\int_{-\left(  2k+1\right)  \pi}^{\left(  2k+1\right)  \pi
}v^{1-p}\left(  x,t\right)  \cos xdx=0\
\]
and thus
\[
\int_{2k\pi}^{\left(  2k+1\right)  \pi}v^{1-p}\left(  x,t\right)  \cos
xdx>0\ \ \ \text{for all\ \ \ }t\in\lbrack0,T_{\max}).
\]
However, by (\ref{decr})$\ $we have$\ $%
\[
\int_{2k\pi}^{\left(  2k+1\right)  \pi}v^{1-p}\left(  x,t\right)  \cos
xdx<0\ \ \ \text{for all\ \ \ }t\in\lbrack0,T_{\max}),
\]
which gives a contradiction.$%
\hfill
\square$

\begin{theorem}
\label{thmC}(\textsf{convergence of type-two blow-up for }$p\geq
2\ $\textsf{with symmetric} $v_{0}\left(  x\right)  $) Let$\ \Phi\left(
x\right)  =\cos x\ $on $\left[  -\pi/2,\pi/2\right]  \ $and $\Phi\left(
x\right)  =0\ $otherwise.\ Assume $v_{0}\left(  x\right)  >0\in C^{\infty
}\left(  S_{m}^{1}\right)  $ satisfies condition\ (\ref{integral-cond})\ and
(\ref{v0}).$\ $Then there exists a sequence of times $t_{n}\nearrow T_{\max}$
such that
\begin{equation}
\lim_{n\rightarrow\infty}\frac{v\left(  x,t_{n}\right)  }{v\left(
0,t_{n}\right)  }=\Phi\left(  x\right)  \ \ \ \text{uniformly on\ \ \ }%
x\in\left[  -m\pi,m\pi\right]  . \label{cos}%
\end{equation}

\end{theorem}

\begin{remark}
Note that if we have type-one blow-up, then\ we consider the rescaling
$v\left(  x,t\right)  /R\left(  t\right)  ,\ $where $R\left(  t\right)
\ $is\ comparable\ to $v_{\max}\left(  t\right)  .$ Hence here for type-two
blow-up, by analogy, it is reasonable to look at the rescaling$\ v\left(
x,t\right)  /v_{\max}\left(  t\right)  ,$ which is (\ref{cos}).
\end{remark}

%

\proof
By Lemma \ref{lem4}, $v\left(  x,t\right)  \ $has type-two blow-up and so
$\left(  T_{\max}-t\right)  ^{1/p}v_{\max}\left(  t\right)  \ $is not
bounded\ on$\ t\in\lbrack0,T_{\max}).\ $Hence there exists a sequence
$s_{n}\nearrow T_{\max}$ such that%
\begin{equation}
\lim_{n\rightarrow\infty}\left(  T_{\max}-s_{n}\right)  ^{1/p}v_{\max}\left(
s_{n}\right)  =\infty. \label{59}%
\end{equation}
Let%
\begin{equation}
\psi_{n}\left(  x\right)  =\frac{1}{T_{\max}-s_{n}}\int_{s_{n}}^{T_{\max}%
}\frac{v\left(  x,s\right)  }{v\left(  0,s\right)  }ds,\ \ \ x\in\left[
-m\pi,m\pi\right]  . \label{60}%
\end{equation}
As we shall be interested in the behavior of $\psi_{n}\left(  x\right)  $ for
$n$ large, without loss of generality, we may assume that $v_{\max}\left(
t\right)  =v\left(  0,t\right)  \ $is increasing in time for all $t\in
\lbrack0,T_{\max})$ (see Lemma \ref{lem2}) and by (\ref{vvv})\ we have
\begin{equation}
0<\psi_{n}\left(  x\right)  \leq1\ \ \ \text{and\ \ }\ \left\vert \psi
_{n}^{\prime}\left(  x\right)  \right\vert \leq1\ \ \ \text{for all\ \ \ }%
x\in\left[  -m\pi,m\pi\right]  \label{psi}%
\end{equation}
for all $n.\ $We also have$\ \psi_{n}\left(  x\right)  =\psi_{n}\left(
-x\right)  \ $for all $x\in\left[  0,m\pi\right]  \ $and $n.\ $Moreover we
have for all $n\ $that
\begin{equation}
\psi_{n}^{\prime}\left(  x\right)  <0\ \ \ \text{for all\ \ \ }x\in\left(
0,m\pi\right)  . \label{PSI}%
\end{equation}
Let $0<K<\pi/2\ $be a fixed number but arbitrary.\ By (\ref{vvcos}) we know
that when $t$ is close to $T_{\max},\ $there holds%
\begin{equation}
v\left(  x,t\right)  \geq v_{\max}\left(  t\right)  \cos x=v\left(
0,t\right)  \cos x\ \ \ \text{for\ all}\ \ \ x\in\left[  -K,K\right]  .
\label{vvvcos}%
\end{equation}
and so$\ $%
\begin{equation}
\lim_{t\rightarrow T_{\max}}v\left(  x,t\right)  =\infty\ \ \ \text{for\ all}%
\ \ \ x\in\left[  -K,K\right]  . \label{vvinf}%
\end{equation}
Moreover, by Lemma \ref{lem2}, we also have%
\[
\left(  v_{xx}+v\right)  \left(  x,t\right)  >0,\ \ \ x\in\left[  -K,K\right]
\]
when $t$ is close to $T_{\max}.\ $As a consequence, when $n$ is large,$\ $we
have%
\begin{align}
0  &  <\frac{1}{T_{\max}-s_{n}}\int_{s_{n}}^{T_{\max}}\frac{v_{xx}\left(
x,s\right)  +v\left(  x,s\right)  }{v\left(  0,s\right)  }ds=\frac{1}{T_{\max
}-s_{n}}\int_{s_{n}}^{T_{\max}}\frac{v_{s}\left(  x,s\right)  }{v^{p}\left(
x,s\right)  v\left(  0,s\right)  }ds\nonumber\\
&  \leq\frac{1}{T_{\max}-s_{n}}\int_{s_{n}}^{T_{\max}}\frac{v_{s}\left(
x,s\right)  }{v^{p+1}\left(  x,s\right)  }ds=\frac{1}{p\left(  T_{\max}%
-s_{n}\right)  v^{p}\left(  x,s_{n}\right)  }\ \ \ \text{for all\ \ \ }%
x\in\left[  -K,K\right]  . \label{psi-2}%
\end{align}

By (\ref{psi}), we may assume that $\psi_{n}\left(  x\right)  $ converges
uniformly on $S_{m}^{1}\ $to a some$\ w\left(  x\right)  \in C^{0}\left(
S_{m}^{1}\right)  \ $and $w\left(  x\right)  \geq0\ $in$\ S_{m}^{1}.$ For any
test function $\varphi\in C_{0}^{\infty}\left(  -\pi/2,\pi/2\right)  $, choose
$0<K<\pi/2\ $so that $\left(  -K,K\right)  \ $contains the support of
$\varphi.\ $By Fubini theorem and integration by parts we have
\begin{align}
&  \int_{-\frac{\pi}{2}}^{\frac{\pi}{2}}w\left(  x\right)  \left[
\varphi_{xx}\left(  x\right)  +\varphi\left(  x\right)  \right]
dx=\lim_{n\rightarrow\infty}\int_{-K}^{K}\psi_{n}\left(  x\right)  \left[
\varphi_{xx}\left(  x\right)  +\varphi\left(  x\right)  \right]  dx\nonumber\\
&  =\lim_{n\rightarrow\infty}\int_{-K}^{K}\left[  \left(  \frac{1}{T_{\max
}-s_{n}}\int_{s_{n}}^{T_{\max}}\frac{v\left(  x,s\right)  }{v\left(
0,s\right)  }ds\right)  \left[  \varphi_{xx}\left(  x\right)  +\varphi\left(
x\right)  \right]  \right]  dx\nonumber\\
&  =\lim_{n\rightarrow\infty}\frac{1}{T_{\max}-s_{n}}\int_{s_{n}}^{T_{\max}%
}\left[  \int_{-K}^{K}\frac{v\left(  x,s\right)  }{v\left(  0,s\right)
}\left[  \varphi_{xx}\left(  x\right)  +\varphi\left(  x\right)  \right]
dx\right]  ds\nonumber\\
&  =\lim_{n\rightarrow\infty}\frac{1}{T_{\max}-s_{n}}\int_{s_{n}}^{T_{\max}%
}\left[  \int_{-K}^{K}\frac{v_{xx}\left(  x,s\right)  +v\left(  x,s\right)
}{v\left(  0,s\right)  }\varphi\left(  x\right)  dx\right]  ds\nonumber\\
&  =\lim_{n\rightarrow\infty}\int_{-K}^{K}\left[  \frac{1}{T_{\max}-s_{n}}%
\int_{s_{n}}^{T_{\max}}\frac{v_{xx}\left(  x,s\right)  +v\left(  x,s\right)
}{v\left(  0,s\right)  }ds\right]  \varphi\left(  x\right)  dx=0
\end{align}
due to (\ref{59}),\ (\ref{vvvcos}) and\ (\ref{psi-2}). This implies that
$w\left(  x\right)  $ is a weak solution of the ODE $w_{xx}+w=0$ in $\left(
-\pi/2,\pi/2\right)  \ $(note that since $\left\vert \psi_{n}^{\prime}\left(
x\right)  \right\vert $ is uniformly bounded, the function $w\ $is Lipschitz
continuous\ with $w\in W^{1,2}\left(  S_{m}^{1}\right)  $))$.$ Regularity
theory implies that $w\left(  x\right)  $ is smooth\ in $x\in\left(
-\pi/2,\pi/2\right)  $ with $w_{xx}+w=0.$ By our definition, $\psi_{n}\left(
x\right)  $ is decreasing in $x$ for $x\in\left(  0,m\pi\right)  $ and has a
maximum at $x=0$ with $\psi_{n}\left(  0\right)  =1.$ This implies that
$w\left(  x\right)  $ is decreasing for $x\in\left(  0,m\pi\right)  $ and has
a maximum at $x=0.\ $Hence $w\left(  0\right)  =1,$ $w^{\prime}\left(
0\right)  =0,$\ and therefore $w\left(  x\right)  =\cos x$ for $x\in\left(
-\pi/2,\pi/2\right)  .$

By Lemma \ref{lem2-1}, we have for large $n$
\begin{align}
0  &  \leq\frac{1}{T_{\max}-s_{n}}\int_{s_{n}}^{T_{\max}}\int_{-K}^{K}\left(
\frac{v\left(  x,s\right)  }{v\left(  0,s\right)  }-\cos x\right)
dxds\nonumber\\
&  =\int_{-K}^{K}\frac{1}{T_{\max}-s_{n}}\int_{s_{n}}^{T_{\max}}\left(
\frac{v\left(  x,s\right)  }{v\left(  0,s\right)  }-\cos x\right)
dsdx\nonumber\\
&  =\int_{-K}^{K}\left(  \psi_{n}\left(  x\right)  -\cos x\right)
dx\rightarrow0\ \ \ \text{as\ \ \ }n\rightarrow\infty. \label{100}%
\end{align}
Hence if we let%
\[
f\left(  s\right)  =\int_{-K}^{K}\left(  \frac{v\left(  x,s\right)  }{v\left(
0,s\right)  }-\cos x\right)  dx,\ \ \ s\in\lbrack0,T_{\max})
\]
we would have
\[
0\leq\frac{1}{T_{\max}-s_{n}}\int_{s_{n}}^{T_{\max}}f\left(  s\right)
ds\rightarrow0\ \ \ \text{as\ \ \ }n\rightarrow\infty.
\]
Therefore by mean value theorem we can find a sequence $s_{n}^{\prime}%
,\ s_{n}<s_{n}^{\prime}<T_{\max},$ so that%
\[
\int_{-K}^{K}\left(  \frac{v\left(  x,s_{n}^{\prime}\right)  }{v\left(
0,s_{n}^{\prime}\right)  }-\cos x\right)  dx\rightarrow0\ \ \ \text{as\ \ \ }%
n\rightarrow\infty.
\]
Note that both $\cos x$ and $v\left(  x,s_{n}^{\prime}\right)  /v\left(
0,s_{n}^{\prime}\right)  $ are bounded functions with bounded derivatives (and
their bounds are independent\ of $n$), and also
\begin{equation}
\frac{v\left(  x,s_{n}^{\prime}\right)  }{v\left(  0,s_{n}^{\prime}\right)
}-\cos x\geq0\ \ \ \text{on\ \ \ }\left[  -K,K\right]  \label{KKK}%
\end{equation}
for large $n.$ Thus\ by Arzela-Ascoli theorem we must have $v\left(
x,s_{n}^{\prime}\right)  /v\left(  0,s_{n}^{\prime}\right)  \rightarrow\cos x$
(passing to a subsequence if necessary)\ uniformly on $\left[  -K,K\right]  $
as $n\rightarrow\infty.\ $

Let $K_{j}$ be a sequence with$\ K_{j}\rightarrow\pi/2$ as $j\rightarrow
\infty.$ For each $j,$ there is a sequence $s_{n}^{\left(  j\right)  }$ so
that $v\left(  x,s_{n}^{\left(  j\right)  }\right)  /v\left(  0,s_{n}^{\left(
j\right)  }\right)  \rightarrow\cos x$ uniformly on $\left[  -K_{j}%
,K_{j}\right]  $ as $j\rightarrow\infty.$ By a diagonal argument, there is a
sequence $\lambda_{n}\nearrow T_{\max}$ such that $v\left(  x,\lambda
_{n}\right)  /v\left(  0,\lambda_{n}\right)  \ $converges uniformly to $\cos
x\ $on $\left[  -K,K\right]  $ for any $0<K<\pi/2.\ $

To obtain the convergence (\ref{cos})\ on $\left[  -\pi/2,\pi/2\right]  $, we
argue as follows (for convenience, any further subsequence of $\lambda_{n}$ is
still denoted as $\lambda_{n}$).\ Assume $v\left(  x,\lambda_{n}\right)
/v\left(  0,\lambda_{n}\right)  \ $does not converge uniformly to $\cos x\ $on
$\left[  -\pi/2,\pi/2\right]  $. Then there exist\ $\varepsilon>0,\ $a
sequence of points $x_{n}\in\left[  -\pi/2,\pi/2\right]  ,$\ and a time
subsequence $\lambda_{n},$ so that
\begin{equation}
f\left(  x_{n},\lambda_{n}\right)  :=\frac{v\left(  x_{n},\lambda_{n}\right)
}{v\left(  0,\lambda_{n}\right)  }-\cos x_{n}\geq\varepsilon\ \ \ \text{for
all\ \ \ }n. \label{f}%
\end{equation}
By the above discussion we may assume that $x_{n}\rightarrow\pi/2.\ $Now by
mean value theorem and (\ref{vvv})%
\begin{align}
\varepsilon &  <\frac{v\left(  x_{n},\lambda_{n}\right)  }{v\left(
0,\lambda_{n}\right)  }\leq\frac{\left\vert v\left(  x_{n},\lambda_{n}\right)
-v\left(  x_{n}-\varepsilon/100,\lambda_{n}\right)  \right\vert }{v\left(
0,\lambda_{n}\right)  }+\frac{v\left(  x_{n}-\varepsilon/100,\lambda
_{n}\right)  }{v\left(  0,\lambda_{n}\right)  }\nonumber\\
&  \leq\frac{\varepsilon}{100}+\frac{v\left(  x_{n}-\varepsilon/100,\lambda
_{n}\right)  }{v\left(  0,\lambda_{n}\right)  } \label{ff}%
\end{align}
where $v\left(  x_{n}-\varepsilon/100,\lambda_{n}\right)  /v\left(
0,\lambda_{n}\right)  \rightarrow\cos\left(  \pi/2-\varepsilon/100\right)
\ $as $n\rightarrow\infty.\ $We have got a contradiction.

Since $v\left(  x,\lambda_{n}\right)  /v\left(  0,\lambda_{n}\right)  \ $is
decreasing in $x\in\left(  0,m\pi\right)  \ $for each time\ $\lambda_{n}$,\ it
must converge to zero uniformly outside the interval $\left[  -\pi
/2,\pi/2\right]  .\ $The proof of Theorem \ref{thmC} is\ done.$%
\hfill
\square\ $

\ \ \ \ \ \ \ \ 

We next want to improve Theorem \ref{thmC} and show that the convergence in
(\ref{cos})\ is valid for all $t\rightarrow T_{\max},\ $not just along a
sequence of times $t_{n}\nearrow T_{\max}.\ $In below, we basically follow
similar arguments as in Lemmas 4.4, 4.5, and 4.6 of Friedman-McLeod\ \cite{FM}%
\ and look more closely at the solution behavior. These estimates are
interesting on their own also.

In the following we still assume that the initial data$\ v_{0}\left(
x\right)  $ satisfies the symmetric condition (\ref{v0}).

\begin{lemma}
\label{lem-new-1}If $x\in\left(  \pi/2,\pi\right)  ,\ $then%
\begin{equation}
\frac{d}{dt}\int_{0}^{x}v^{1-p}\left(  y,t\right)  \cos ydy<0\ \ \ \text{for
all\ \ \ }t\in\lbrack0,T_{\max}). \label{new1}%
\end{equation}

\end{lemma}%

\proof
We proceed as in \cite{FM}, Lemma 4.5.\ By direct computation, we have%
\[
\frac{d}{dt}\int_{0}^{x}v^{1-p}\left(  y,t\right)  \cos ydy=\left(
1-p\right)  \left[  v_{x}\left(  x,t\right)  \cos x+v\left(  x,t\right)  \sin
x\right]  <0
\]
since for $x\in\left(  \pi/2,\pi\right)  $ we have $v_{x}\left(  x,t\right)
<0,\ \cos x<0,\ \sin x>0.%
\hfill
\square$

\begin{lemma}
\label{lem-new-2}If $x>\pi/2\ $or $x<-\pi/2,$ then there exists a constant $C$
depending on $x$ such that%
\begin{equation}
0<v\left(  x,t\right)  \leq C\ \ \ \text{for all\ \ \ }t\in\lbrack0,T_{\max})
\label{new2}%
\end{equation}
i.e.,\ $v\left(  x,t\right)  $ does not blow up for $\left\vert x\right\vert
>\pi/2.$
\end{lemma}

%

\proof
We proceed as in \cite{FM}, Lemma 4.6.\ Since$\ v\left(  x,t\right)  $ is
decreasing in $x\in\left(  0,m\pi\right)  \ $for all time,$\ $without loss of
generality, we may just look at the case $x\in\left(  \pi/2,\pi\right)
.\ $Suppose$\ v\left(  x,t\right)  \ $is not bounded, then there exists a
sequence $t_{n}\nearrow T_{\max}$ so that $v\left(  x,t_{n}\right)
\rightarrow\infty.$ By Lemma \ref{lem2} we must have $v\left(  x,t\right)
\rightarrow\infty$ as $t\rightarrow T_{\max}.\ $In particular, we have (note
that $v\left(  y,t\right)  $ is decreasing for $y>0$)%
\begin{equation}
\int_{0}^{x}v^{1-p}\left(  y,t\right)  \cos ydy\rightarrow
0\ \ \ \text{as\ \ \ }t\rightarrow T_{\max}. \label{new3}%
\end{equation}
On the other hand, we may write for fixed small $\delta>0$%
\begin{align}
&  \int_{0}^{x}v^{1-p}\left(  y,t\right)  \cos ydy\nonumber\\
&  =\int_{0}^{\left(  \pi-\delta\right)  /2}v^{1-p}\left(  y,t\right)  \cos
ydy+\int_{\left(  \pi-\delta\right)  /2}^{\left(  \pi+\delta\right)
/2}v^{1-p}\left(  y,t\right)  \cos ydy+\int_{\left(  \pi+\delta\right)
/2}^{x}v^{1-p}\left(  y,t\right)  \cos ydy. \label{new4}%
\end{align}
As$\ v\left(  x,t\right)  >0$ is decreasing in $x\in\left(  0,m\pi\right)  ,$
the second term in (\ref{new4})\ is negative for all time$.\ $Also by
(\ref{vvcos})\ in Lemma \ref{lem2-1}, there is a constant $c>0\ $such that
$v\left(  y,t\right)  \geq cv\left(  0,t\right)  $ for all $y\in\left[
0,\left(  \pi-\delta\right)  /2\right]  \ $and all time large enough.\ Finally
for $y\in\left[  \left(  \pi+\delta\right)  /2,x\right]  ,\ $by Theorem
\ref{thmC}\ there exists a\ sequence $t_{n}\nearrow T_{\max}\ $such that
$v\left(  y,t_{n}\right)  \leq\varepsilon_{n}v\left(  0,t_{n}\right)  ,$ where
$\varepsilon_{n}\rightarrow0\ $as $n\rightarrow\infty.\ $Hence$\ $for\ $n$
large enough we conclude
\begin{align*}
&  \int_{0}^{x}v^{1-p}\left(  y,t_{n}\right)  \cos ydy\\
&  \leq c^{1-p}v^{1-p}\left(  0,t_{n}\right)  \int_{0}^{\left(  \pi
-\delta\right)  /2}\cos ydy+\varepsilon_{n}^{1-p}v^{1-p}\left(  0,t_{n}%
\right)  \int_{\left(  \pi+\delta\right)  /2}^{x}\cos ydy<0.
\end{align*}
This gives a contradiction due to (\ref{new1}) and (\ref{new3}).

The proof for the case $x<-\pi/2$ is similar. $%
\hfill
\square$

\begin{lemma}
\label{lem-new-3}Let $w\left(  x\right)  $ be a nonnegative\ Lipschitz
function defined on $\left[  -\pi/2,\pi/2\right]  .\ $Suppose that
$w\ $satisfies the inequality $w_{xx}+w\geq0\ $on$\ \left(  -\pi
/2,\pi/2\right)  $ in the sense of distribution. If $w\left(  -\pi/2\right)
=w\left(  \pi/2\right)  =0,$ then $w\left(  x\right)  =a\cos x,\ $where
$a=\max_{x\in\left[  -\pi/2,\pi/2\right]  }w\left(  x\right)  .$
\end{lemma}

%

\proof
Let $\varphi_{n}$ be a sequence of smooth nonnegative functions with compact
support in $\left(  -\pi/2,\pi/2\right)  \ $such that it converges to $w$ in
$H_{0}^{1}\left[  -\pi/2,\pi/2\right]  \ $(note that$\ 0\leq w\in H_{0}%
^{1}\left[  -\pi/2,\pi/2\right]  $)$.\ $Since $w\ $satisfies $w_{xx}+w\geq
0\ $on$\ \left(  -\pi/2,\pi/2\right)  $ in the sense of distribution, we have%
\[
\int_{-\frac{\pi}{2}}^{\frac{\pi}{2}}\left(  \varphi_{n}^{\prime}\left(
x\right)  w^{\prime}\left(  x\right)  -\varphi_{n}\left(  x\right)  w\left(
x\right)  \right)  dx\leq0\ \ \ \text{for all\ \ \ }n.
\]
Letting $n\rightarrow\infty$ we get%
\begin{equation}
\int_{-\frac{\pi}{2}}^{\frac{\pi}{2}}\left(  \left(  w^{\prime}\left(
x\right)  \right)  ^{2}-w^{2}\left(  x\right)  \right)  dx\leq0. \label{ww1}%
\end{equation}
Note that $\lambda=1$ is the principal eigenvalue of the operator
$d^{2}/dx^{2}$ on the interval$\ \left[  -\pi/2,\pi/2\right]  \ $with $\cos x$
the principal eigenfunction satisfying Dirichlet boundary condition. Thus we
obtain%
\begin{equation}
\int_{-\frac{\pi}{2}}^{\frac{\pi}{2}}w^{2}\left(  x\right)  dx\leq\int
_{-\frac{\pi}{2}}^{\frac{\pi}{2}}\left(  w^{\prime}\left(  x\right)  \right)
^{2}dx \label{ww2}%
\end{equation}
where equality holds only when $w$ is a constant multiple of the principal
eigenfunction.\ Equations (\ref{ww1}) and (\ref{ww2})\ imply that $w$ is a
principal eigenfunction on the interval $\left[  -\pi/2,\pi/2\right]  $ and so
$w_{xx}+w=0\ $on$\ \left(  -\pi/2,\pi/2\right)  .\ $Let $a=\max_{x\in\left[
-\pi/2,\pi/2\right]  }w\left(  x\right)  .$ Then we conclude that $w=a\cos
x\ $on$\ \left[  -\pi/2,\pi/2\right]  .%
\hfill
\square$

\begin{theorem}
\label{thmC-1}Under the same assumption as in Theorem \ref{thmC} we have%
\begin{equation}
\lim_{t\rightarrow T_{\max}}\frac{v\left(  x,t\right)  }{v\left(  0,t\right)
}=\Phi\left(  x\right)  \ \ \ \text{uniformly on\ \ \ }x\in\left[  -m\pi
,m\pi\right]  . \label{cosine}%
\end{equation}

\end{theorem}%

\proof
It suffices to prove that for any sequence $t_{j}\nearrow T_{\max}$ there is a
subsequence, also denoted as $t_{j},$ so that%
\[
\lim_{j\rightarrow\infty}\frac{v\left(  x,t_{j}\right)  }{v\left(
0,t_{j}\right)  }=\Phi\left(  x\right)  \ \ \ \text{uniformly on\ \ \ }%
x\in\left[  -m\pi,m\pi\right]  .
\]
This would imply that the convergence is for all time $t\rightarrow T_{\max
}.\ $Let $t_{j}$ be a sequence with $t_{j}\nearrow T_{\max}.$ By Lemma
\ref{lem2} there is a subsequence $t_{j}$ and a nonnegative
Lipschitz\ function $w\left(  x\right)  $ defined on $\left[  -m\pi
,m\pi\right]  $ so that
\[
\lim_{j\rightarrow\infty}\frac{v\left(  x,t_{j}\right)  }{v\left(
0,t_{j}\right)  }=w\left(  x\right)  \ \ \ \text{uniformly on\ \ \ }%
x\in\left[  -m\pi,m\pi\right]  .
\]
We clearly have $\max_{x\in\left[  -\pi/2,\pi/2\right]  }w=1\ $and by Lemma
\ref{lem2} it satisfies $w_{xx}+w\geq0\ $on$\ \left(  -\pi/2,\pi/2\right)  $
in the sense of distribution. By Lemma \ref{lem-new-2}, since $v\left(
x,t\right)  $ does not blow up for $\left\vert x\right\vert >\pi/2,\ $we must
have $w\left(  x\right)  =0\ $for $\left\vert x\right\vert >\pi/2.$ By
continuity, we have$\ w\left(  -\pi/2\right)  =w\left(  \pi/2\right)  =0.$
Thus Lemma \ref{lem-new-3} implies that $w\left(  x\right)  =\cos x\ $for
$\left\vert x\right\vert <\pi/2.\ $The proof is done.$%
\hfill
\square$

\subsection{Convergence to a translational self-similar solution.}

Back to the slow-speed curve contracting flow $\left(  \bigstar\right)  $
with$\ \alpha\in(0,1]$, the initial curve $\gamma_{0}$ has curvature
$k_{0}\left(  x\right)  >0$ satisfying (\ref{integral-cond}), i.e.,%
\begin{equation}
\int_{S_{m}^{1}}\frac{1}{k_{0}\left(  x\right)  }e^{ix}dx=0.
\end{equation}
If $k_{0}\left(  x\right)  $ satisfies the symmetric condition
\begin{equation}
k_{0}\left(  x\right)  =k_{0}\left(  -x\right)  \ \ \ \text{and\ \ \ }%
k_{0}^{\prime}\left(  x\right)  <0,\text{\ \ \ \ \ }x\in\left(  0,m\pi\right)
\end{equation}
and $k_{\max}\left(  t\right)  \ $of$\ \gamma_{t}$ has type-two blow-up, then
by Theorem \ref{thmC-1} we have the convergence
\begin{equation}
\lim_{t\rightarrow T_{\max}}\frac{k\left(  x,t\right)  }{k\left(  0,t\right)
}=\left(  \cos x\right)  ^{\frac{1}{\alpha}}\ \ \ \text{uniformly
on\ \ \ }x\in\left[  -\pi/2,\pi/2\right]  .
\end{equation}

It is well known that under the flow $\left(  \bigstar\right)  $, there is a
special \textbf{translational self-similar solution} $\Gamma_{t}\ $translating
in the direction $\left(  0,1\right)  $ with unit speed\ (see \cite{NT} or the
book\ \cite{CZ}). For each time $t,\ \Gamma_{t}$ is only a translation of
$\Gamma_{0}\ $(this $\Gamma_{0}$ is not a closed curve, but still convex and
embedded). If we use tangent angle $x$ to parametrize $\Gamma_{0},$ its
parametrization$\ $is given by%
\[
\Gamma_{0}=\left(  \int_{0}^{x}\frac{\cos\xi}{\left(  \cos\xi\right)
^{\frac{1}{\alpha}}}d\xi,\ \int_{0}^{x}\frac{\sin\xi}{\left(  \cos\xi\right)
^{\frac{1}{\alpha}}}d\xi\right)  ,\ \ \ x\in\left(  -\pi/2,\pi/2\right)
\]
where
\[
\int_{0}^{x}\frac{\sin\xi}{\left(  \cos\xi\right)  ^{\frac{1}{\alpha}}}%
d\xi=\left\{
\begin{array}
[c]{l}%
\frac{\alpha}{\alpha-1}\left[  1-\left(  \cos x\right)  ^{1-\frac{1}{\alpha}%
}\right]  ,\ \ \ \alpha\in\left(  0,1\right)
\vspace{3mm}%
\\
-\log\cos x,\ \ \ \alpha=1.
\end{array}
\right.
\]
In particular the curve $\Gamma_{0}\ $goes to infinity as $x\rightarrow\pm
\pi/2.$ The curvature of $\Gamma_{0}\ $at angle $x\ $is given by $k\left(
x\right)  =\left(  \cos x\right)  ^{1/\alpha},\ x\in\left(  -\pi
/2,\pi/2\right)  ,$ with maximum at $x=0.\ $When $\alpha=1,$ we get\ Grayson's
"\textbf{Grim Reaper}",\textsf{\ }which is$\ \Gamma_{0}=\left(  x,-\log\cos
x\right)  ,\ x\in\left(  -\pi/2,\pi/2\right)  .$

Evolve the above given symmetric $\gamma_{0}$ according to the flow $\left(
\bigstar\right)  $. For any $t\in\lbrack0,T_{\max}),$ choose the point
$x_{t}\in\gamma_{t}$ at which the curvature is $k_{\max}\left(  t\right)  $
(by the assumption there is only one such point)\ and translate $\gamma_{t}$
so that $x_{t}$ becomes the origin $O=\left(  0,0\right)  $. Call this
translational curve $\tilde{\gamma}_{t}$. Next rotate it so that the unit
tangent vector at the origin of $\tilde{\gamma}_{t}\ $becomes $\left(
1,0\right)  ,$ and finally dilate the curve so that its maximal curvature
becomes $1\ $and denote this final curve as $\hat{\gamma}_{t}.$ Theorem
\ref{thmC-1} says that if we have type-two blow-up of $k_{\max}\left(
t\right)  ,$ then over the region $x\in\left(  -\pi/2,\pi/2\right)
,\ \hat{\gamma}_{t}$ converges to the above translational self-similar
solution\ $\Gamma_{0}$ as $t\rightarrow T_{\max}.$ When $\alpha=1,$ this
phenomenon has been observed by Angenent in \cite{ANG}.\ 

Thus we can summarize the following important observation of the slow speed
flow $\left(  \bigstar\right)  $:\ \textbf{for type-one blow-up, the
asymptotic behavior is given by a homothetic self-similar solution, while for
type-two blow-up, the asymptotic behavior (in the special symmetric case)\ is
given by a translational self-similar solution.}

To end this paper we point out that most of the lemmas and theorems remain
valid even the initial condition $v_{0}\left(  x\right)  $ does not satisfy
the integral condition (\ref{integral-cond}), as long as it is positive,
smooth, and $2m\pi$-periodic.\ They include Lemmas \ref{lem1},\ \ref{lem2}%
,\ \ref{lem2-1},\ \ref{lem-new-1} and\ Theorems \ref{thmB}%
,\ \ref{thm-type-one}, \ref{thmA},

As for Lemma \ref{lem-new-2}\ and Theorems\ \ref{thmC}, \ref{thmC-1}, if we
add the extra assumption that $\left(  T_{\max}-t\right)  ^{1/p}v_{\max
}\left(  t\right)  \ $is not bounded\ on$\ t\in\lbrack0,T_{\max}),\ $then they
are all valid even if $v_{0}\left(  x\right)  $ does not satisfy
(\ref{integral-cond}).

\emph{In particular, we emphasize\ again that for }$p\in\lbrack2,\infty
),\ $\emph{there is either type-one blow-up\ or type-two blow-up.\ Moreover,
type-one blow-up occurs only when}$\ v_{0}\left(  x\right)  \ $%
\emph{satisfies\ the integral condition (\ref{integral-cond}) and if }%
$v_{0}\left(  x\right)  $\emph{ does not satisfy (\ref{integral-cond}), then
the blow-up is always of type-two.\ Thus the generic blow-up behavior for
}$p\in\lbrack2,\infty)$\emph{ is type-two. }

\subsection{What to do next\ ?}

There is still a difficult question of estimating the type-two blow-up rate of
$v\left(  x,t\right)  =k^{\alpha}\left(  x,t\right)  .\ $When $\alpha=1$
(i.e., $p=1+1/\alpha=2$)\ and $m=2,\ $Angenent$\ $%
and\ Vel\'{a}zquez\ \cite{AV} had given a nontrivial proof of the existence of
some symmetric\ initial\ data $v_{0}\left(  x\right)  >0,$ satisfying
(\ref{integral-cond}),\ with the type-two blow-up rate$\ $
\[
v_{\max}\left(  t\right)  =\left(  1+o\left(  1\right)  \right)  \sqrt
{\frac{\ln\ln\left(  \frac{1}{T_{\max}-t}\right)  }{T_{\max}-t}}%
\ \ \ \text{as\ \ \ }t\rightarrow T_{\max}%
\]
and therefore
\begin{equation}
v_{\max}\left(  t\right)  \sqrt{T_{\max}-t}\sim\sqrt{\ln\ln\left(  \frac
{1}{T_{\max}-t}\right)  }\rightarrow\infty\;\;\;\text{as\ \ \ }t\rightarrow
T_{\max}. \label{Ve}%
\end{equation}
We are wondering if certain similar estimate holds in the case $\alpha
\in(0,1]\ $(i.e., $p>2$).\ At this moment we do not know and we hope to work
on it in the future.\ 

\section{Some pictures for the ODE\ (\ref{wode})}

In this section we give some pictures relating to the ODE$\ w^{\prime\prime
}+w-w^{1-p}=0.\ $These pictures can help us understand convergence behavior
(for general $p\in\left(  -\infty,\infty\right)  $) of the PDE $\partial
u/\partial\tau=u^{p}\left(  u_{xx}+u-u^{1-p}\right)  \ $(with positive initial
data$\ u_{0}\in C^{\infty}\left(  S_{m}^{1}\right)  $ and periodic boundary
condition).\ This is because that the ODE is a steady state of the PDE.
Let$\ $%
\[
F\left(  s\right)  =\left\{
\begin{array}
[c]{l}%
s^{2}-\frac{2}{2-p}s^{2-p},\ \ \ p\neq2,\ \ \ p\in\left(  -\infty
,\infty\right)
\vspace{3mm}%
\\
s^{2}-2\log s,\ \ \ p=2
\end{array}
\right.  ,\ \ \ s\in\left(  0,\infty\right)  .
\]
The graphs of $F\left(  s\right)  $ for $p=-1\in\left(  -\infty,0\right)
,\ p=1\in\left(  0,2\right)  ,\ p=3\in\lbrack2,\infty)$ are given below:%
\[%
{\includegraphics[
height=1.8282in,
width=2.4284in
]%
{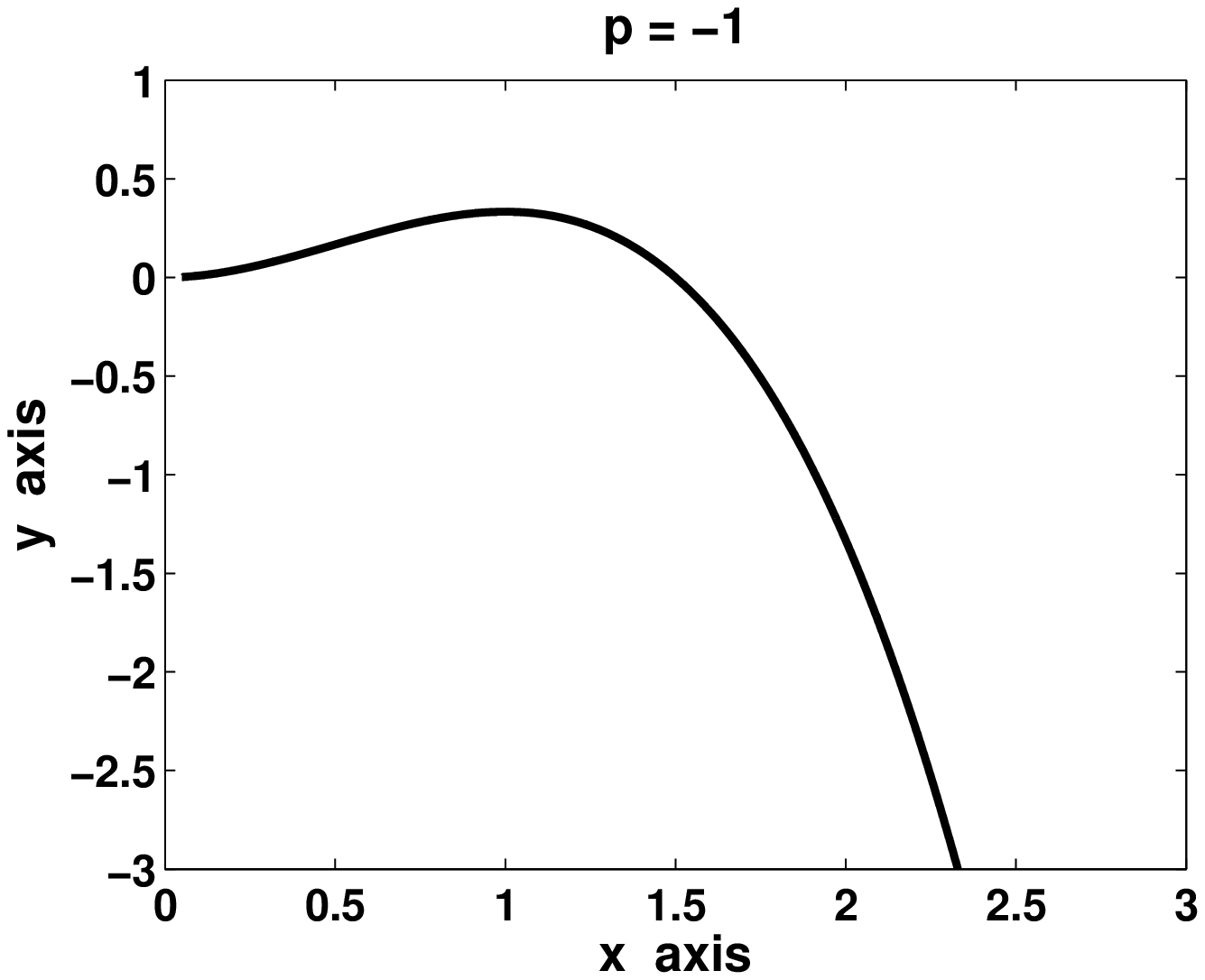}%
}%
\]

\begin{center}

$%
\begin{array}
[c]{cl}%
{\includegraphics[
height=1.8862in,
width=2.5054in
]%
{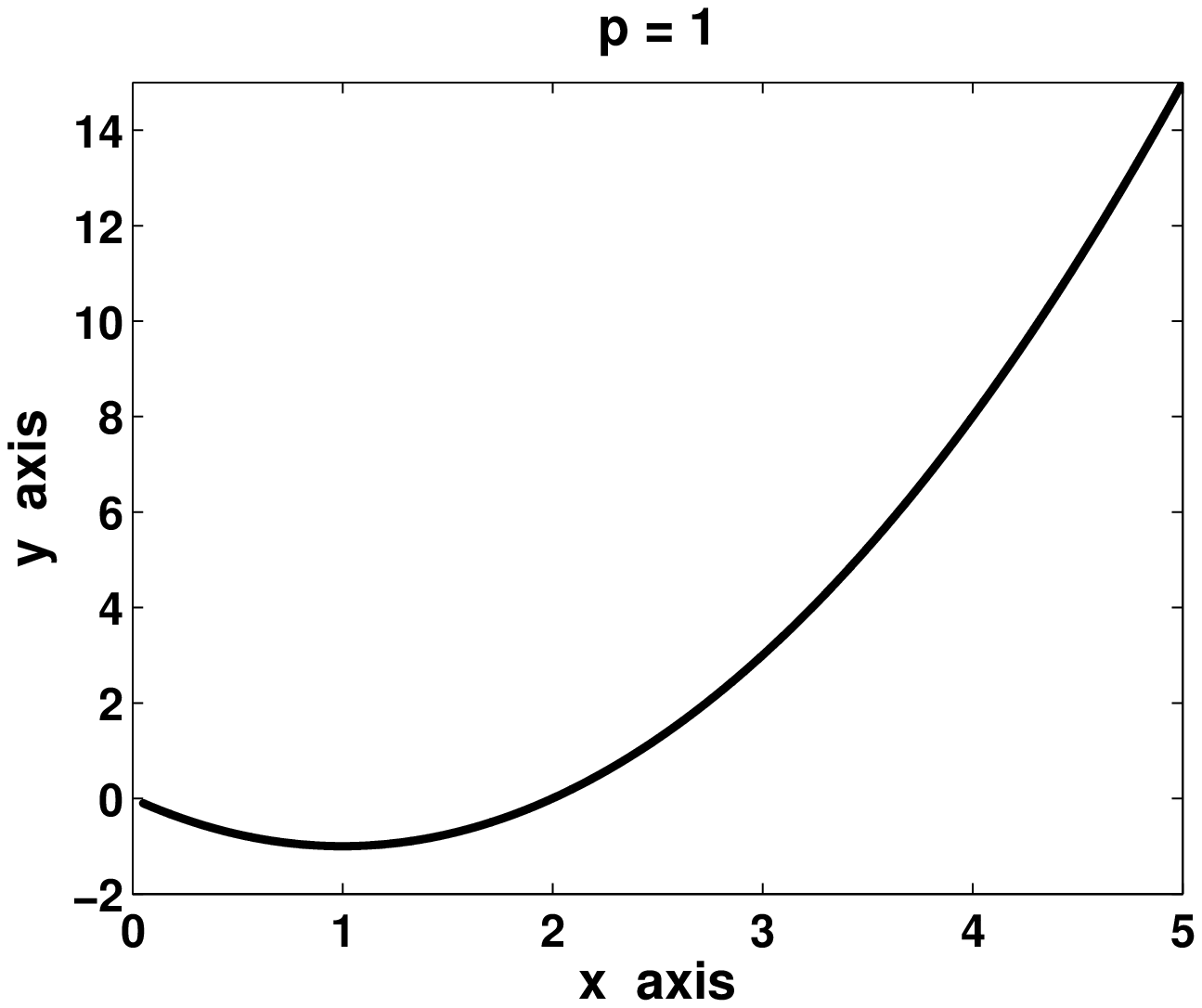}%
}%
&
{\includegraphics[
height=1.881in,
width=2.4984in
]%
{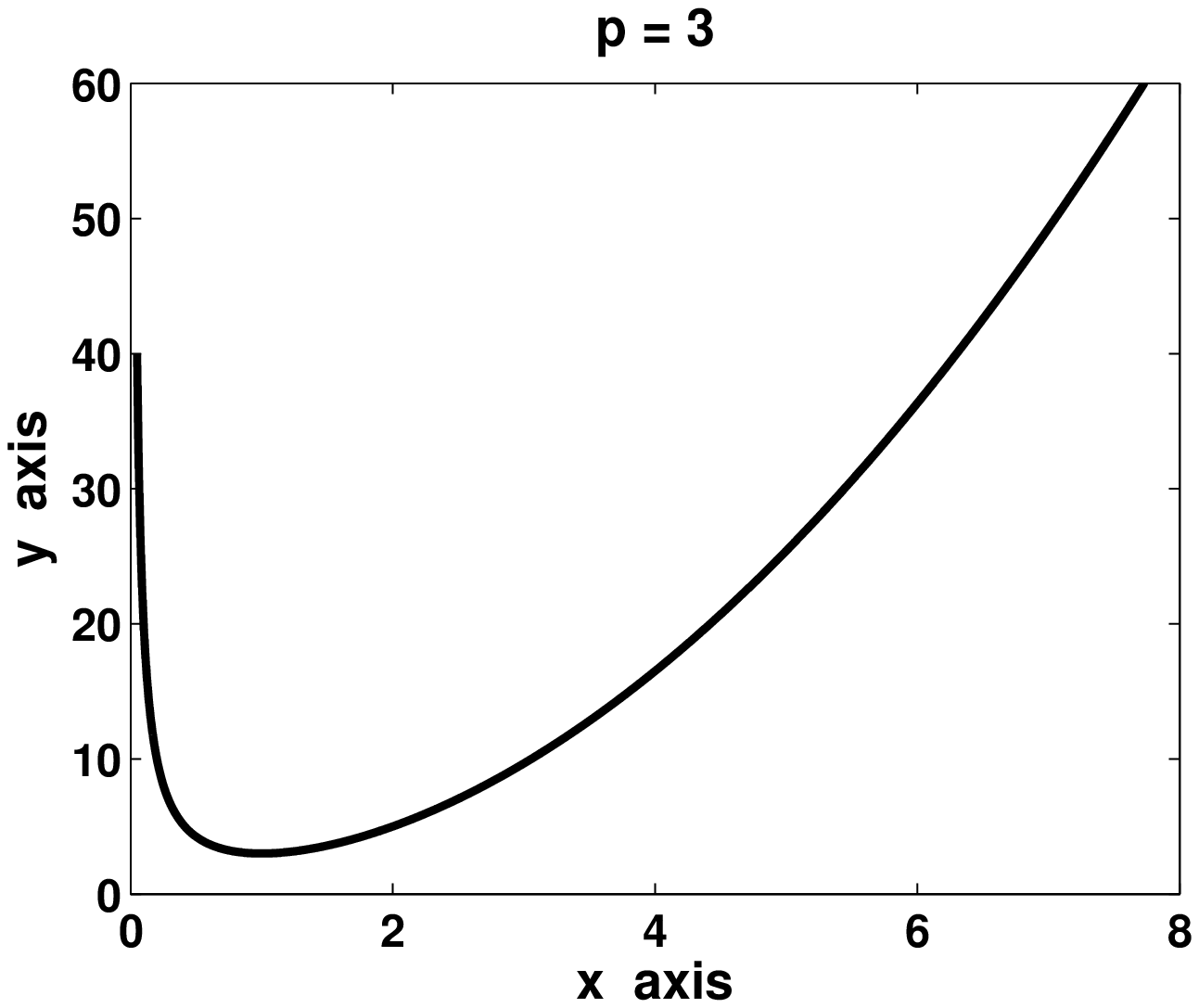}%
}%
\end{array}
\ $
\end{center}

The graph of $F\left(  s\right)  $ for any $p\in\left(  -\infty,0\right)
\ $is analogous to the above picture for $p=-1.\ F\left(  s\right)  $
increases on $s\in\left(  0,1\right)  \ $with$\ F\left(  0\right)  =0,$
$F\left(  1\right)  >0,$\ and decreases to $-\infty\ $on $\left(
1,\infty\right)  .\ $The graph of $F\left(  s\right)  $ for any $p\in\left(
0,2\right)  \ $is analogous to the above picture for $p=1.\ F\left(  s\right)
$ decreases on $s\in\left(  0,1\right)  \ $with$\ F\left(  0\right)
=0,\ F\left(  1\right)  <0,$\ and increases to $+\infty\ $on $\left(
1,\infty\right)  .\ $Finally the graph of $F\left(  s\right)  $ for any
$p\in\lbrack2,\infty)\ $is analogous to the above picture for $p=3.\ F\left(
s\right)  $ decreases on $s\in\left(  0,1\right)  \ $with$\ \lim
_{s\rightarrow0^{+}}F\left(  s\right)  =+\infty,$ $F\left(  1\right)
>0,$\ and increases to $+\infty\ $on $\left(  1,\infty\right)  $. Also note
that when $p=0,$ $F\left(  s\right)  \equiv0.\ $

Any solution $w\left(  x\right)  \ $to the ODE $w^{\prime\prime}+w-w^{1-p}=0$
must satisfy the identity$\ w_{x}^{2}\left(  x\right)  =F\left(  M\right)
-F\left(  w\left(  x\right)  \right)  \ $for all $x\ $in the domain $I\ $of
$w\left(  x\right)  ,\ $on which $w\left(  x\right)  >0.\ $Here we may assume
$0\in I\ $and $w\left(  0\right)  =M\geq1$\ is the maximum value of $w$ on
$I.\ $For $p\in\left(  -\infty,0\right)  ,$ we only have type-one blow up for
$v\left(  x,t\right)  $ of equation $\left(  \clubsuit\right)  $.\ If the
rescaled$\ $solution$\ u\left(  x,t\right)  =v\left(  x,t\right)  /R\left(
t\right)  $ converges to $w\left(  x\right)  $ on some interval $I,\ $then we
must have $w\left(  x\right)  \equiv1\ $over $I$.$\ $Otherwise we have$\ M>1$
and use the first picture to get%
\begin{equation}
w_{x}^{2}\left(  x\right)  =F\left(  M\right)  -F\left(  w\left(  x\right)
\right)  <0
\end{equation}
for all\ $x\in I\ $such that $1\leq w\left(  x\right)  <M.\ $This gives a
contradiction and so $w\left(  x\right)  \equiv1\ $over $I.\ $

The main difference between $p\in\left(  0,2\right)  \ $and $p\in
\lbrack2,\infty)\ $is that there exist bump solutions (degenerate)\ to the
ODE\ for$\ p\in\left(  0,2\right)  ,$ but for $p\in\lbrack2,\infty),\ $all
solutions to the ODE are positive everywhere\ and periodic over $\mathbb{R}$
(nondegenerate). Again, this can also be seen from the second and third pictures.

\ \ \ \ 

\bigskip

\textbf{Acknowledgments.\ \ \ }While writing this paper, we had discussions
with several mathematicians including Professors Ben Andrews,\ Sigurd
Angenent,\ Hiroshi Matano, Jong-Sheng Guo and Chia-Hsing Nien.\ We are very
grateful to all of them.\ The third author would like to acknowledge the
support of the National Science Council and the National Center for
Theoretical Sciences\ of Taiwan.

\bigskip

\bigskip\ 

\bigskip\ \ \

\ \ \ \ 

\ 

\ \ \ \ 

Chi-Cheung Poon

Department of Mathematics, National Chung Cheng University, Chiayi 621,\ TAIWAN.

Email:\ \textit{ccpoon@math.ntu.edu.tw}

\ \ \ \ \ \ 

\ \ \ \ \ 

\ \ \ \ 

Yu-Chu Lin and\ Dong-Ho Tsai

Department of Mathematics, National Tsing Hua University, Hsinchu 300,\ TAIWAN.

Email:\ \textit{yclin@math.nthu.edu.tw,\ \ \ dhtsai@math.nthu.edu.tw}

\end{document}